\documentclass[twoside,leqno,twocolumn]{article}
\usepackage{ltexpprt}
\usepackage{amsfonts,amssymb}
\usepackage{color}
\usepackage{pifont}
\usepackage{multirow}
\usepackage{amsmath}
\usepackage{url}
\usepackage{graphicx}
 \usepackage{mathrsfs}
\usepackage{cite}
\usepackage{enumitem}
 \usepackage{cases}
 \usepackage{pbox}
\usepackage{setspace}

 \usepackage[usenames,dvipsnames]{pstricks}
 \usepackage{epsfig}
 \usepackage{pst-grad} 
 \usepackage{pst-plot} 

\usepackage{makecell}
\usepackage[margin=1in]{geometry}
\usepackage{relsize}
\usepackage{fancyhdr}

\newcommand{\bcap} {\hspace{2pt} \mathlarger{\cap}
\hspace{2pt}}

    \makeatletter
    \def\@listi{\leftmargin\leftmargini
        \parsep 1\p@ \@plus0\p@ \@minus\p@
        \topsep 2\p@   \@plus0\p@ \@minus\p@
        \itemsep1\p@ \@plus0\p@ \@minus\p@}
    \let\@listI\@listi\@listi
    \makeatother

\newcommand*\mcapinn[2]{\vcenter{\hbox{$\mathsurround=0pt
\ifx\displaystyle#1\textstyle\else#1\fi\bigcap$}}}

\newcommand*\mcupinn[2]{\vcenter{\hbox{$
\bigcup$}}}

\newcommand{\bP}[1]{{\mathbb{P}}\left[{#1}\right]}

\begin{document}

\title{On $k$-Connectivity and Minimum Vertex Degree\\in Random $s$-Intersection
Graphs}
\author{Jun Zhao \hspace{49pt} Osman
Ya\u{g}an \hspace{49pt} Virgil Gligor\vspace{3pt}\\
Carnegie Mellon University\thanks{The authors are with ECE
Department and CyLab, Carnegie Mellon University, USA. Emails:
\{junzhao, oyagan, virgil\}@andrew.cmu.edu}\vspace{-5pt}}
\date{}

\maketitle

\thispagestyle{fancy}

\fancyhead[C]{{ACM-SIAM Meeting on Analytic Algorithmics and Combinatorics (ANALCO) 2015}}

\fancyhead[L]{}
\fancyhead[R]{}

\pagestyle{fancy}


\begin{abstract} {\small
\emph{Random $s$-intersection graphs} have recently received much interest
in a wide range of application areas. Broadly speaking, a random
$s$-intersection graph is constructed by first assigning each vertex
a set of items in some \emph{random} manner, and then putting an undirected
edge between all pairs of vertices that share at least $s$ items
(the graph is called a \emph{random intersection graph} when $s=1$). A
special case of particular interest is a \emph{uniform random
$s$-intersection graph}, where each vertex independently selects the
same number of items uniformly at random from a common item pool.
Another important case is a \emph{binomial random $s$-intersection graph},
where each item from a pool is independently assigned to each vertex
with the same probability. Both models have found numerous
applications thus far including cryptanalysis, and the modeling of
recommender systems, secure sensor networks, online social networks,
trust networks and small-world networks (uniform random
$s$-intersection graphs), as well as clustering analysis,
classification, and the design of integrated circuits (binomial
random $s$-intersection graphs).

In this paper, for binomial/uniform random $s$-intersection graphs,
we present results related to $k$-connectivity and minimum vertex
degree. Specifically, we derive the asymptotically exact
probabilities and zero--one laws for the following three properties:
(i) $k$-vertex-connectivity, (ii) $k$-edge-connectivity and (iii)
the property of minimum vertex degree being at least
$k$.}\vspace{3pt}
\end{abstract}

\vspace{-3pt}

\section{Introduction}

\emph{Random $s$-intersection graphs} have received considerable
attention recently
\cite{Rybarczyk,bloznelis2013,Assortativity,mil10,ZhaoYaganGligor,Perfectmatchings,JZISIT14,Bloznelis201494,ball2014,fullver,r4,Y,X}.
In such a graph, each vertex is equipped with a set of items in some \emph{random}
manner, and two vertices establish an undirected edge in between if
and only if they have at least $s$ items in common. A large amount
of work
\cite{r1,herdingRKG,PES:6114960,YaganThesis,Shang,ZhaoAllerton,virgil,GodehardtJaworski,ryb3,zz,2013arXiv1301.0466R,CohenThesis,ZhaoCDC,DBLP:journals/corr/abs-1301-7320,Jaworski20062152,ISIT,BloznelisD13,Nikoletseas:2008:LIS:1414105.1414429,RSA:RSA20005,Karonski99}
study the case of $s$ being $1$, under which the graphs are simply
referred to as \emph{random intersection graphs}.

Random ($s$-)intersection graphs have been used to model secure
wireless sensor networks
\cite{Rybarczyk,JZISIT14,adrian,ISIT,virgil,YaganThesis,qcomp_kcon,ZhaoAllerton,ANewell,zz,2013arXiv1301.0466R},
wireless frequency hopping
 \cite{ZhaoAllerton}, epidemics in human populations \cite{ball2014,mil10}, small-world networks \cite{5383986}, trust networks \cite{virgillncs,Ysb}, social networks \cite{Assortativity,bloznelis2013,mil10,ZhaoYaganGligor,r4} such as collaboration networks \cite{Assortativity,bloznelis2013,mil10} and common-interest networks \cite{ZhaoYaganGligor,r4}.
Random intersection graphs also motivated Beer \emph{et al.}
\cite{beer2011vertex-journal,beer2011vertex-conf} to introduce a
general concept of {\em vertex random graphs} that subsumes any
graph model where {\em random} features are assigned to vertices,
and edges are drawn based on deterministic relations between the
features of the vertices.

 Among different models of random $s$-intersection graphs, two widely studied models are the so-called \emph{uniform random $s$-intersection graph} and \emph{binomial random $s$-intersection graph} defined in detail below.

\subsection{Graph models}

\paragraph{Uniform random $s$-intersection graph.}
A \emph{uniform random $s$-intersection graph}, denoted by
$G_s(n,K_n,P_n)$, is defined on $n$ vertices as follows. Each vertex
\emph{independently} selects $K_n$ different items \emph{uniformly
at random} from a pool of $P_n$ distinct items. Two vertices have an
edge in between
 if and only if they share at least $s$ items.
 The notion ``uniform'' means that all vertices have the same number of items (but likely different sets of items). Here $K_n$ and $P_n$ are both functions of $n$, while $s$ does not scale with $n$. It holds that $1\leq s\leq K_n \leq
 P_n$. Under $s=1$, the graph is also known as a
 \emph{random key graph} \cite{mobihocQ1,yagan,5383986}.

\paragraph{Binomial random $s$-intersection graph.}
A \emph{binomial random $s$-intersection graph}, denoted by
$H_s(n,t_n,P_n)$, is defined on $n$ vertices as follows. Each item
from a pool of $P_n$ distinct items is assigned to each vertex
\emph{independently} with probability $t_n$. Two vertices establish
an edge in between
 if and only if they have at least $s$ items in
common. The term ``binomial'' is used since the number of items
assigned to each vertex follows a binomial distribution with
parameters $P_n$ (the number of trials) and $t_n$ (the success
probability in each trial). Here $t_n$ and $P_n$ are both functions
of $n$, while $s$ does not scale with $n$. Also it holds that $1\leq
s \leq P_n$.

\subsection{Problem Statement.}

Our goal in this paper is to investigate  properties related to
\emph{$k$-connectivity} and \emph{minimum vertex degree} of random
$s$-intersection graphs (\emph{$k$-vertex-connectivity} and \emph{$k$-edge-connectivity} are called together as
 $k$-connectivity for convenience).
In particular, we wish to answer the following question:

For a uniform random $s$-intersection graph
 $G_s(n,K_n,P_n)$ (resp., a binomial random $s$-intersection graph $H_s(n,t_n,P_n)$), with parameters $K_n$ (resp., $t_n$) and $P_n$ scaling with the number of vertices $n$,
what is the asymptotic behavior of the probabilities for
$G_s(n,K_n,P_n)$ (resp., $H_s(n,t_n,P_n)$) (i) being
$k$-vertex-connected, (ii) being $k$-edge-connected, and (iii) having a
minimum vertex degree at least $k$, respectively, as $n$ grows
large? A graph is said to be $k$-vertex-connected if the remaining
graph is still connected despite the deletion of at most $(k-1)$
arbitrary vertices, and $k$-edge-connectivity is defined similarly
for the deletion of edges \cite{Bollobas}; with $k=1$, these
definitions reduce to the standard notion of {\em graph
connectivity} \cite{citeulike:4012374,Z}. The degree of a vertex is
defined as the number of edges incident on it. The three graph
properties considered here are related to each other in that
$k$-vertex-connectivity implies $k$-edge-connectivity, which in turn
implies that the minimum vertex degree is at least $k$
\cite{Bollobas}.

\subsection{Summary of Results.}

We summarize our results below, first for a uniform random
$s$-intersection graph and then for a binomial random
$s$-intersection graph. Throughout the paper, both $s$ and $k$ are
positive integers and do not scale with $n$. Also, naturally we
consider $1 \leq s \leq K_n \leq P_n$ for graph $G_s(n,K_n,P_n)$ and
$1 \leq s \leq P_n$ for graph $H_s(n,t_n,P_n)$. We use the standard
Landau asymptotic notation $\Omega(\cdot), \omega(\cdot), O(\cdot),
o(\cdot),\Theta(\cdot)$. $\mathbb{P}[\mathcal {E}]$ denotes the
probability that event $\mathcal {E}$ happens.

\textbf{$k$-Connectivity \& minimum vertex degree in uniform random
$s$-intersection graphs:}
\\
For a uniform random $s$-intersection graph $G_s(n,K_n,P_n)$ under
$P_n = \Omega(n)$, with sequence $\alpha_n$ defined by
\begin{align}
\frac{1}{s!} \cdot \frac{{K_n}^{2s}}{{P_n}^{s}}  &  = \frac{\ln  n +
{(k-1)} \ln \ln n + {\alpha_n}}{n}, \label{crit}
\end{align}
then as $n \to \infty$, if $ \alpha_n \to \alpha^{\star} \in
[-\infty, \infty]$, the following convergence results hold:\\$
\mathbb{P} \left[\hspace{1.5pt}G_s(n,K_n,P_n)\textrm{ is
$k$-vertex-connected}.\hspace{1pt}\right] \to e^{-
\frac{e^{-\alpha^{\star}}}{(k-1)!}}$,\\
$ \mathbb{P} \left[\hspace{1pt}G_s(n,K_n,P_n) \textrm{ is
$k$-edge-connected}.\hspace{1.5pt}\right] \to e^{-
\frac{e^{-\alpha^{\star}}}{(k-1)!}}$,\\and \vspace{-3pt}
\begin{align}
  \mathbb{P}\left[
  \begin{array}{l} \hspace{-3pt} G_s(n,K_n,P_n) \textrm{ has a minimum\hspace{-3pt}}
  \\  \hspace{-3pt}\textrm{vertex degree at least }k.
\end{array} \right] & \to e^{- \frac{e^{-\alpha^{\star}}}{(k-1)!}} . ~~~~~~~~~~~~~~~\nonumber
\end{align}

\textbf{$k$-Connectivity \& minimum vertex degree in binomial random
$s$-intersection graphs:}
\\
For a binomial random $s$-intersection graph $H_s(n,t_n,P_n)$ under
$P_n = \Omega(n)$ for $s\geq 2$ or $P_n = \Omega(n^c)$ for $s=1$
with some constant $c>1$, with sequence $\beta_n$ defined by
\begin{align}
\frac{1}{s!} \cdot {t_n}^{2s}{P_n}^{s}  &  = \frac{\ln  n + {(k-1)}
\ln \ln n + {\beta_n}}{n}, \nonumber
\end{align}
then as $n \to \infty$, if $ \beta_n \to \beta^{\star} \in [-\infty,
\infty]$, the following convergence results hold:
\\$ \mathbb{P}
\left[\hspace{1.5pt}H_s(n,t_n,P_n)\textrm{ is
$k$-vertex-connected}.\hspace{1.5pt}\right] \to
e^{- \frac{e^{-\beta^{\star}}}{(k-1)!}}$, \\
$ \mathbb{P} \left[\hspace{1.5pt}H_s(n,t_n,P_n) \textrm{ is
$k$-edge-connected}.\hspace{1.5pt}\right] \to e^{-
\frac{e^{-\beta^{\star}}}{(k-1)!}}$,\\and \vspace{-3pt}
\begin{align}
  \mathbb{P}\left[
  \begin{array}{l} \hspace{-3pt} H_s(n,t_n,P_n) \textrm{ has a minimum\hspace{-3pt}}
  \\  \hspace{-3pt}\textrm{vertex degree at least }k.
\end{array} \right] & \to e^{- \frac{e^{-\beta^{\star}}}{(k-1)!}} . ~~~~~~~~~~~~~~~\nonumber
\end{align}

Since the probability $e^{- \frac{e^{-\alpha^{\star}}}{(k-1)!}}$
(resp., $e^{- \frac{e^{-\beta^{\star}}}{(k-1)!}}$) equals $1$ if
$\alpha^{\star} = \infty$ (resp., $\beta^{\star} = \infty$) and $0$
if $ \alpha^{\star} = -\infty$ (resp., $ \beta^{\star} = -\infty$),
the above results of asymptotically exact probabilities also imply
the corresponding zero--one laws, where a zero--one law \cite{yagan}
means that the probability that the graph has  certain property
asymptotically converges to $0$ under some conditions and converges
to $1$ under some other conditions.

\subsection{Comparison with related work.}

Table \ref{table:related-work} summarizes relevant work in the
literature on uniform/binomial random $s$-intersection graphs in
terms of $k$-vertex-connectivity, $k$-edge connectivity, and the
property of minimum vertex degree being at least $k$.

Among the related work, Bloznelis and Rybarczyk
\cite{Bloznelis201494} recently also derived the asymptotically exact
probabilities of \emph{uniform} random $s$-intersection graphs (but not of \emph{binomial} random $s$-intersection graphs) for the three properties above (the easily implied results on $k$-edge-connectivity were not explicitly mentioned). Yet, when $s$ is a constant or $O(1)$ as in many
applications, their results require $K_n = O\big((\ln
n)^{\frac{1}{5s}}\big)$ for $k$-connectivity ({$k$-vertex-connectivity} and {$k$-edge-connectivity}), under a scaling the same as in
Equation (\ref{crit}). In other words, the 
one-law part of $k$-connectivity is as follows: under certain conditions including $K_n = O\big((\ln
n)^{\frac{1}{5s}}\big)$, if $\frac{1}{s!} \cdot
\frac{{K_n}^{2s}}{{P_n}^{s}}  = \frac{\ln  n + {(k-1)} \ln \ln n +
\omega(1)}{n}$, then a uniform random $s$-intersection graph
$G_s(n,K_n,P_n)$ is $k$-connected with a probability converging to $1$
as $n\to\infty$. From $K_n = O\big((\ln n)^{\frac{1}{5s}}\big)$ and
$\frac{1}{s!} \cdot \frac{{K_n}^{2s}}{{P_n}^{s}}  = \frac{\ln  n +
{(k-1)} \ln \ln n + \omega(1)}{n}$, it is straightforward to derive
$P_n = O\big(n^{\frac{1}{s}} (\ln n)^{-\frac{3}{5s}}\big)$; i.e.,
$P_n =\widetilde{O}\big(n^{\frac{1}{s}}\big)$ ignoring the $\ln n$
terms.
%
%
%
%
 However, in secure wireless sensor
network applications where uniform random $s$-intersection graphs
are widely investigated, conditions $K_n = O\big((\ln
n)^{\frac{1}{5s}}\big)$ and $P_n
=\widetilde{O}\big(n^{\frac{1}{s}}\big)$ are both likely impractical
because $K_n$ and $P_n$ are often at least on the order of $\ln n$
and $n$, respectively, to ensure that the network has reasonable
resiliency against sensor capture attacks
\cite{adrian,virgil,YaganThesis}. The
results reported in this paper cover the practical range where $P_n$
is at least on the order of $n$.

\begin{table*}[t]
\vspace{2mm} \normalsize  \vspace{3pt}
\begin{center}
\begin{tabular}{!{\vrule width 1.15pt}l|l|l|l|l!{\vrule width 1.15pt}}
   \Xhline{2.5\arrayrulewidth}
\multicolumn{2}{!{\vrule width 1.15pt}l|}{\hspace{-1.5pt}Graph}
& \hspace{-1.5pt}Property
& \hspace{-1.5pt}Results                            &
\hspace{-1.5pt}Work
\\ \hline
\multirow{9}{*}{\begin{tabular}[c]{@{}l@{}}\hspace{-1.5pt}uniform
random\\ \hspace{-1.5pt}$s$-intersection
graph\hspace{-1.5pt}\end{tabular}}  &
\multirow{3}{*}{\hspace{-1.5pt}$G_s(n,K_n,P_n)$\hspace{-1.5pt}} &
\multirow{2}{*}{\begin{tabular}[c]{@{}l@{}}\hspace{-1.5pt}$k$-connectivity
\& \\ \hspace{-1.5pt}min. vertex degree $\geq k$\end{tabular}} &
\multirow{2}{*}{\hspace{-1.5pt}exact probabilities\hspace{-1.5pt}} &
\hspace{-1.5pt}\textbf{this paper}
\\ \cline{5-5}
                                                                                              &                                   &                                                                                                             &                                    & \begin{tabular}[c]{@{}l@{}}\hspace{-1.5pt}\cite{Bloznelis201494}~(only for\\ \hspace{-1.5pt}$K_n = O\big((\ln n)^{\frac{1}{5s}}\big)$)\end{tabular}  \\ \cline{3-5}
                                                                                              &                                   & \begin{tabular}[c]{@{}l@{}}\hspace{-1.5pt}connectivity \& \\ \hspace{-1.5pt}min. vertex degree $\geq 1$\end{tabular}                      & \hspace{-1.5pt}exact probabilities                  & \begin{tabular}[c]{@{}l@{}}\hspace{-1.5pt}\cite{Perfectmatchings}~(only for\\  \hspace{-1.5pt}$K_n = O\big((\ln n)^{\frac{1}{5s}}\big)$)\end{tabular} \\ \cline{2-5}
                                                                                              & \multirow{4}{*}{\hspace{-1.5pt}$G_1(n,K_n,P_n)$\hspace{-1.5pt}} & \multirow{2}{*}{\begin{tabular}[c]{@{}l@{}}\hspace{-1.5pt}$k$-connectivity \& \\ \hspace{-1.5pt}min. vertex degree $\geq k$\hspace{-1.5pt}\end{tabular}} & \hspace{-1.5pt}exact probabilities\hspace{-1.5pt}                  & \hspace{-1.5pt}\cite{ZhaoCDC}                                                                                  \\ \cline{4-5}
                                                                                              &                                   &                                                                                                             & \hspace{-1.5pt}zero--one laws                       & \hspace{-1.5pt}\cite{zz,ISIT}                                                                             \\ \cline{3-5}
                                                                                              &                                   & \multirow{2}{*}{\begin{tabular}[c]{@{}l@{}}\hspace{-1.5pt}connectivity \& \\ \hspace{-1.5pt}min. vertex degree $\geq 1$\hspace{-1.5pt}\end{tabular}}     & \hspace{-1.5pt}exact probabilities\hspace{-1.5pt}                  & \hspace{-1.5pt}\cite{ryb3}                                                                                     \\ \cline{4-5}
                                                                                              &                                   &                                                                                                             & \hspace{-1.5pt}zero--one laws                       & \hspace{-1.5pt}\cite{r1,yagan}                                                                                 \\ \hline
\multirow{8}{*}{\begin{tabular}[c]{@{}l@{}}\hspace{-1.5pt}binomial
random\\ \hspace{-1.5pt}$s$-intersection
graph\hspace{-1.5pt}\end{tabular}} &
\multirow{4}{*}{\hspace{-1.5pt}$H_s(n,t_n,P_n)$\hspace{-1.5pt}} &
\multirow{2}{*}{\begin{tabular}[c]{@{}l@{}}\hspace{-1.5pt}$k$-connectivity
\& \\ \hspace{-1.5pt}min. vertex degree $\geq
k$\hspace{-1.5pt}\end{tabular}} & \hspace{-1.5pt}exact
probabilities\hspace{-1.5pt}                  &
\multirow{4}{*}{\hspace{-1.5pt}\textbf{this paper}\hspace{-1.5pt}}
\\ \cline{4-4}
                                                                                              &                                   &                                                                                                             & \hspace{-1.5pt}zero--one laws                       &                                                                                                   \\ \cline{3-4}
                                                                                              &                                   & \multirow{2}{*}{\begin{tabular}[c]{@{}l@{}}\hspace{-1.5pt}connectivity \& \\ min. vertex degree $\geq 1$\hspace{-1.5pt}\end{tabular}}     & \hspace{-1.5pt}exact probabilities\hspace{-1.5pt}                  &                                                                                                   \\ \cline{4-4}
                                                                                              &                                   &                                                                                                             & \hspace{-1.5pt}zero--one laws                       &                                                                                                   \\ \cline{2-5}
                                                                                              & \multirow{4}{*}{\hspace{-1.5pt}$H_1(n,t_n,P_n)$\hspace{-1.5pt}} & \multirow{2}{*}{\begin{tabular}[c]{@{}l@{}}\hspace{-1.5pt}$k$-connectivity \& \\ \hspace{-1.5pt}min. vertex degree $\geq k$\hspace{-1.5pt}\end{tabular}} & \hspace{-1.5pt}exact probabilities\hspace{-1.5pt}                   & \hspace{-1.5pt}\cite{ZhaoCDC}                                                                                                  \\ \cline{4-5}
                                                                                              &                                   &                                                                                                             & \hspace{-1.5pt}zero--one laws\hspace{-1.5pt}                       & \hspace{-1.5pt}\cite{zz}                                                                                       \\ \cline{3-5}
                                                                                              &                                   & \multirow{2}{*}{\begin{tabular}[c]{@{}l@{}}\hspace{-1.5pt}connectivity \& \\ \hspace{-1.5pt}min. vertex degree $\geq 1$\hspace{-1.5pt}\end{tabular}}     & \hspace{-1.5pt}exact probabilities\hspace{-1.5pt}                  & \hspace{-1.5pt}\cite{2013arXiv1301.0466R}                                                                      \\ \cline{4-5}
                                                                                              &                                   &                                                                                                             & \hspace{-1.5pt}zero--one laws                       & \hspace{-1.5pt}\cite{CohenThesis,Shang}                                                                           \\    \Xhline{2.5\arrayrulewidth}
\end{tabular}\vspace{2pt}
 \caption{Comparison of our results with related work. $k$-vertex connectivity and $k$-edge connectivity are together written as
 $k$-connectivity. Note that results on $k$-connectivity and min. vertex degree $\geq k$ also
 imply the corresponding results on connectivity and min. vertex degree $\geq 1$ by setting $k$ as $1$.
 Also, results of exact probabilities imply the corresponding results of zero--one laws
 by monotonicity arguments.\vspace{-5pt}} \label{table:related-work}
  \end{center}
\end{table*}


\subsection{Roadmap.}

We organize the rest of the paper as follows. We detail the main
results in Section \ref{sec:res}. Sections \ref{sec:basic:ideas} and
\ref{sec:prf:thm:bin} detail the steps of establishing the theorems.
We conclude the paper in Section \ref{sec:Conclusion}. The Appendix
provides additional arguments used in proving the theorems.

\section{Main Results} \label{sec:res}

Below we explain the main results of uniform random $s$-intersection
graphs and binomial random $s$-intersection graphs, respectively.

\subsection{Results of uniform random $s$-intersection graphs.}

The following theorem presents results on $k$-connectivity and
minimum vertex degree in a uniform random $s$-intersection graph
$G_s(n,K_n,P_n)$.

\begin{theorem} \label{thm:uni}
For a uniform random $s$-intersection graph $G_s(n,K_n,P_n)$ under
\vspace{-3pt}
\begin{align}
P_n = \Omega(n),\vspace{-3pt} \label{eqPnOmegan}
\end{align}
with sequence $\alpha_n$ defined by
\begin{align}
\frac{1}{s!} \cdot \frac{{K_n}^{2s}}{{P_n}^{s}}  &  = \frac{\ln  n +
{(k-1)} \ln \ln n + {\alpha_n}}{n},
\vspace{-5pt}\label{thm:uni:eq:edge}
\end{align}
it holds that \vspace{-4pt}
\begin{align}
&  \lim\limits_{n \to \infty} \mathbb{P}[\hspace{2pt} G_s(n,K_n,P_n)
\textrm{ is
  $k$-vertex-connected.} \hspace{2pt} ]
\nonumber \\ &   =  \lim\limits_{n \to \infty}
\mathbb{P}[\hspace{2pt} G_s(n,K_n,P_n) \textrm{ is
  $k$-edge-connected.} \hspace{2pt} ]
\nonumber \\ &   = \lim\limits_{n \to \infty} \mathbb{P}\left[
  \begin{array}{l} \hspace{-4pt} G_s(n,K_n,P_n) \textrm{ has a minimum\hspace{-4pt}}
  \\  \hspace{-4pt}\textrm{vertex degree at least }k.
\end{array} \right] \nonumber \\ &   = \begin{cases} e^{- \frac{e^{-\alpha ^*}}{(k-1)!}}, \hspace{0.1cm} \textrm{if
}\lim\limits_{n \to \infty}{\alpha_n} = \alpha^* \in
(-\infty,\infty),
 \\0, \hspace{1.1cm} \textrm{if }\lim\limits_{n \to
\infty}{\alpha_n} = -\infty,   \\ 1, \hspace{1.1cm} \textrm{if
}\lim\limits_{n \to \infty}{\alpha_n} = \infty.
\end{cases} \nonumber
\end{align}
\end{theorem}

For graph $G_s(n,K_n,P_n)$, Theorem \ref{thm:uni} presents the
asymptotically exact probabilities and zero--one laws for the
following three properties: (i) $k$-vertex connectivity, (ii)
$k$-edge-connectivity and (iii) the property of minimum vertex
degree being at least $k$.

By Lemma \ref{qn-dist} on Page \pageref{qn-dist}, under
(\ref{eqPnOmegan}) and (\ref{thm:uni:eq:edge}) with constrained
$|\alpha_n |= O(\ln \ln n )$, we can show that the left hand side of
(\ref{thm:uni:eq:edge}), i.e., $\frac{1}{s!} \cdot
\frac{{K_n}^{2s}}{{P_n}^{s}} $, is asymptotically equivalent to the
edge probability of graph $G_s(n,K_n,P_n)$. As given in Lemma 13 of
the full version \cite{fullver}, with $q_n$ denoting the edge
probability of graph $G_s(n,K_n,P_n)$, if condition
(\ref{thm:uni:eq:edge}) is replaced by $q_n = \frac{\ln  n + {(k-1)}
\ln \ln n + {\alpha_n}}{n}$,
 and  condition (\ref{eqPnOmegan}) is kept unchanged, then all results in Theorem \ref{thm:uni} still follow. Hence, the uniform random $s$-intersection graph model under condition (\ref{eqPnOmegan})
 exhibits the same behavior as the well-known
 Erd\H{o}s-R\'enyi graph model
  \cite{citeulike:4012374}, in the sense that for each of (i) $k$-vertex-connectivity, (ii) $k$-edge-connectivity and (iii) the property of minimum vertex degree being at least $k$, a common point for the phase transition from a zero-law to a one-law
occurs when the edge probability equals $\frac{\ln  n + {(k-1)} \ln
\ln n}{n}$.

\subsection{Results for binomial random $s$-intersection graphs.}

The following theorem presents results on $k$-connectivity and
minimum vertex degree in a binomial random $s$-intersection graph
$H_s(n,t_n,P_n)$.

\begin{theorem}  \label{thm:bin}
For a binomial random $s$-intersection graph $H_s(n,t_n,P_n)$ under
\begin{align}
\begin{cases}
P_n = \Omega(n), &\textrm{ for } s \geq 2, \\
P_n = \Omega(n^c)\textrm{ for some constant }c>1, &\textrm{ for } s
= 1,
 \end{cases}
\label{thm:bin:eq:P}
\end{align}
with sequence $\beta_n$ defined by
\begin{align}
\frac{1}{s!} \cdot {t_n}^{2s}{P_n}^{s}  &  = \frac{\ln  n + {(k-1)}
\ln \ln n + {\beta_n}}{n},  \label{thm:bin:eq:edge}
\end{align}
it holds that
\begin{align}
&  \lim\limits_{n \to \infty} \mathbb{P}[\hspace{2pt} H_s(n,t_n,P_n)
\textrm{ is
  $k$-vertex-connected.} \hspace{2pt} ]
\nonumber \\ &   =  \lim\limits_{n \to \infty}
\mathbb{P}[\hspace{2pt} H_s(n,t_n,P_n) \textrm{ is
  $k$-edge-connected.} \hspace{2pt} ]
\nonumber \\ &   = \lim\limits_{n \to \infty} \mathbb{P}\left[
  \begin{array}{l} \hspace{-4pt} H_s(n,t_n,P_n) \textrm{ has a minimum\hspace{-4pt}}
  \\  \hspace{-4pt}\textrm{vertex degree at least }k.
\end{array} \right] \nonumber \\ &   = \begin{cases} e^{- \frac{e^{-\beta ^*}}{(k-1)!}}, \hspace{0.1cm} \textrm{if
}\lim\limits_{n \to \infty}{\beta_n} = \beta^* \in (-\infty,\infty),
 \\0, \hspace{1.1cm} \textrm{if }\lim\limits_{n \to
\infty}{\beta_n} = -\infty,   \\ 1, \hspace{1.1cm} \textrm{if
}\lim\limits_{n \to \infty}{\beta_n} = \infty.
\end{cases} \nonumber
\end{align}
\end{theorem}

For graph $H_s(n,t_n,P_n)$, Theorem \ref{thm:bin} presents the
asymptotically exact probabilities and zero--one laws for the
following three properties: (i) $k$-vertex connectivity, (ii)
$k$-edge-connectivity and (iii) the property of minimum vertex
degree being at least $k$.

By Lemma  12 of the full version \cite{fullver}, under
(\ref{thm:bin:eq:P}) and (\ref{thm:bin:eq:edge}) with constrained
$|\beta_n| = O(\ln \ln n  )$, we can show that the left hand side of
(\ref{thm:bin:eq:edge}), i.e., $\frac{1}{s!} \cdot
{t_n}^{2s}{P_n}^{s}  $, is asymptotically equivalent to the edge
probability of graph $H_s(n,t_n,P_n)$. As given in Lemma 14 of the
full version \cite{fullver}, with $\rho_n$ denoting the edge
probability of graph $H_s(n,t_n,P_n)$, if condition
(\ref{thm:bin:eq:edge}) is replaced by $\rho_n = \frac{\ln  n +
{(k-1)} \ln \ln n + {\beta_n}}{n}$,
 and  condition (\ref{thm:bin:eq:P}) is kept unchanged, then all results in Theorem \ref{thm:bin} still follow.
 Therefore, the binomial random $s$-intersection graph model under condition (\ref{thm:bin:eq:P}) exhibits the same behavior with Erd\H{o}s-R\'enyi graph model, in the sense that for each of (i) $k$-vertex-connectivity, (ii) $k$-edge-connectivity, and (iii) the property of minimum vertex degree being at least $k$, a common point for the phase transition from a zero-law to a one-law
occurs when the edge probability equals $\frac{\ln  n + {(k-1)} \ln
\ln n}{n}$.

The condition (\ref{thm:bin:eq:P}) has $P_n = \Omega(n)$ for $s \geq
2$, and requires a stronger one for $s=1$: $P_n = \Omega(n^c)$ for
some constant $c>1$. The range $P_n = \Theta(n)$ is covered by $P_n
= \Omega(n)$, but not by $P_n = \Omega(n^c)$ with $c>1$. For $s=1$
and $P_n = \Theta(n)$, results for $k$-vertex-connectivity, $k$-edge
connectivity, and the property of minimum vertex degree being at
least $k$ use a scaling different from (\ref{thm:bin:eq:edge}), as
given by \cite[Theorem 4]{zz}.

\section{Establishing Theorem \ref{thm:uni}} \label{sec:basic:ideas}

Theorem \ref{thm:uni} in the special case of $s = 1$ is proved by us
\cite{ZhaoCDC}. Below we explain the steps of establishing Theorem
\ref{thm:uni} for $s \geq 2$. In Section \ref{sec:confine:alpha}, we
show that $|\alpha_n|$ can be confined as $O(\ln \ln n) $ in proving
Theorem \ref{thm:uni}. In Section \ref{sec:del:vem}, we consider the
relationships between vertex connectivity, edge connectivity, and
minimum vertex degree.

\subsection{Confining $|\alpha_n|$.} \label{sec:confine:alpha}

To confine $|\alpha_n|$ as $O(\ln \ln n) $ in proving Theorem
\ref{thm:uni}, we will demonstrate
\begin{align}
\textrm{Theorem \ref{thm:uni} under }|\alpha_n | = O(\ln \ln n)
 \Rightarrow \textrm{Theorem \ref{thm:uni}}. \label{cp_alph}
\end{align}
Note that $k$-vertex-connectivity, $k$-edge-connectivity, and the
property of minimum vertex degree being at least $k$ are all
monotone increasing\footnote{A graph property is called monotone
increasing if it holds under the addition of edges
\cite{Bollobas,JansonLuczakRucinski7}.}. For any monotone increasing
property $\mathcal {I}$, the probability that a spanning subgraph
(resp., supergraph) of graph $G$ has $\mathcal {I}$ is at most
(resp., at least) the probability of $G$ having $\mathcal {I}$.
Therefore, to show (\ref{cp_alph}), it suffices to prove the
following lemma.

\begin{lemma} \label{graph_Gs_cpl}

\textbf{(a)} For graph $G_s(n,K_n,P_n)$ under $P_n = \Omega(n)$ and
\begin{align}
\frac{1}{s!} \cdot \frac{{K_n}^{2s}}{{P_n}^{s}}  &  = \frac{\ln  n +
{(k-1)} \ln \ln n + {\alpha_n}}{n} \label{al1-parta}
\end{align}
with $\lim_{n \to \infty}\alpha_n = -\infty$, there exists graph
$G_s(n,\widetilde{K_n},\widetilde{P_n})$ under $\widetilde{P_n} =
\Omega(n)$ and
\begin{align}
\frac{1}{s!} \cdot
\frac{{\widetilde{K_n}}^{2s}}{{\widetilde{P_n}}^{s}}  &  = \frac{\ln
n + {(k-1)} \ln \ln n + {\widetilde{\alpha_n}}}{n} \label{al0-parta}
\end{align}
with $\lim_{n \to \infty}\widetilde{\alpha_n} = -\infty$ and
$\widetilde{\alpha_n} = -O(\ln \ln n)$, such that there exists a
graph coupling\footnote{As used by Rybarczyk
\cite{zz,2013arXiv1301.0466R}, a coupling of two random graphs $G_1$
and $G_2$ means a probability space on which random graphs $G_1'$
and $G_2'$ are defined such that $G_1'$ and $G_2'$ have the same
distributions as $G_1$ and $G_2$, respectively. If $G_1'$ is a
spanning subgraph (resp., supergraph) of $G_2'$, we say that under
the coupling, $G_1$ is a spanning subgraph (resp., supergraph) of
$G_2$, which yields that for any monotone increasing property
$\mathcal {I}$, the probability of $G_1$ having $\mathcal {I}$ is at
most (resp., at least) the probability of $G_2$ having $\mathcal
{I}$.} under which $G_s(n,K_n,P_n)$ is a spanning subgraph of
$G_s(n,\widetilde{K_n},\widetilde{P_n})$.

\textbf{(b)} For graph $G_s(n,K_n,P_n)$ under $P_n = \Omega(n)$ and
\begin{align}
\frac{1}{s!} \cdot \frac{{K_n}^{2s}}{{P_n}^{s}}  &  = \frac{\ln  n +
{(k-1)} \ln \ln n + {\alpha_n}}{n} \label{al1}
\end{align}
with $\lim_{n \to \infty}\alpha_n = \infty$, there exists graph
$G_s(n,\widehat{K_n},\widehat{P_n})$ under $\widehat{P_n} =
\Omega(n)$ and
\begin{align}
\frac{1}{s!} \cdot \frac{{\widehat{K_n}}^{2s}}{{\widehat{P_n}}^{s}}
&  = \frac{\ln  n + {(k-1)} \ln \ln n + {\widehat{\alpha_n}}}{n}
\label{al0}
\end{align}
with $\lim_{n \to \infty}\widehat{\alpha_n} = \infty$ and
$\widehat{\alpha_n} = O(\ln \ln n)$, such that there exists a graph
coupling under which $G_s(n,K_n,P_n)$ is a spanning supergraph of
$G_s(n,\widehat{K_n},\widehat{P_n})$.

\end{lemma}

The proof of Lemma \ref{graph_Gs_cpl} is provided in Section
\ref{sec_graph_Gs_cpl} in the Appendix.

\subsection{Relationships between vertex connectivity, edge connectivity, and minimum vertex degree.} \label{sec:del:vem}

Recall that the vertex connectivity of a graph is defined as the
minimum number of vertices needing to be deleted to have the
remaining graph disconnected, and the edge connectivity is defined
similarly for the deletion of edges \cite{Bollobas}. For graph
$G_s(n,K_n,P_n)$, we use $\kappa_v, \kappa_e$ and $\delta$ to denote
the vertex connectivity, the edge connectivity, and the minimum
vertex degree, respectively. Then $k$-vertex-connectivity,
$k$-edge-connectivity, and
 the property of minimum vertex degree being at least $k$, are given by events $\kappa_v \geq k$, $\kappa_e \geq k$, and $\delta \geq k$, respectively. For any graph, the vertex connectivity is at most the edge connectivity, and the  edge connectivity is at most the minimum vertex degree \cite{Bollobas,ZhaoISIT2014,FJYGISIT2014}. Therefore, $\kappa_v \leq \kappa_e \leq \delta$ holds. Then
\begin{align}
 \mathbb{P}[\hspace{2pt} \kappa_v \geq k\hspace{2pt}] & \leq  \mathbb{P}[ \hspace{2pt}\kappa_e \geq k\hspace{2pt}]  \leq  \mathbb{P}[\hspace{2pt}\delta \geq k\hspace{2pt}],   \label{eq_kpv1}
 \end{align}
and
\begin{align}
 \mathbb{P}[ \kappa_v \geq  k ] &=
  \mathbb{P}[\delta \geq  k ] -
\mathbb{P}[\hspace{2pt}(\kappa_v < k) \cap (\delta \geq k)
\hspace{2pt}] \nonumber \\ &  \geq   \mathbb{P}[\hspace{2pt}\delta
\geq k \hspace{2pt}] - \sum_{\ell = 0}^{k-1}
\mathbb{P}[\hspace{2pt}(\kappa_v = \ell) \cap (\delta > \ell
)\hspace{2pt}].   \label{eq_kpv2}
\end{align}

Therefore, the proof is completed once we show Lemmas \ref{lem:mvd}
and \ref{lem:mvd:con} below. Note that since $k$ is a constant,
condition (\ref{thm:uni:eq:edge}) with $|\alpha_n| = O(\ln \ln n)$
in Theorem \ref{thm:uni} implies condition (\ref{pe-lnn-lnlnn-pm})
in Lemma \ref{lem:mvd:con}.

\begin{lemma}[\hspace{-.1pt}{\cite[Theorem 2]{qcomp_kcon} (our work)}] \label{lem:mvd}

For uniform random $s$-intersection graph $G_s(n,K_n,P_n)$ under
$P_n = \Omega(n)$, if there exists sequence $\alpha_n$ satisfying
$|\alpha_n| = O(\ln \ln n)$ such that
\begin{align}
\frac{1}{s!} \cdot \frac{{K_n}^{2s}}{{P_n}^{s}}  &  = \frac{\ln  n +
{(k-1)} \ln \ln n + {\alpha_n}}{n},\nonumber
\end{align}
then with $\delta$ denoting the minimum vertex degree, it holds that
\begin{align}
 \lim\limits_{n \to \infty}  \mathbb{P}[\delta \geq k] & =
  e^{- \frac{e^{-\alpha^{\star}}}{(k-1)!}}, \textrm{ if }\lim_{n \to \infty} \alpha_n = \alpha^{\star} \in [-\infty, \infty].
 \label{mnd}
\end{align}
\end{lemma}

\begin{lemma} \label{lem:mvd:con}
For uniform random $s$-intersection graph $G_s(n,K_n,P_n)$ under
$P_n = \Omega(n)$ and
\begin{align}
\frac{1}{s!} \cdot \frac{{K_n}^{2s}}{{P_n}^{s}}   &  = \frac{\ln  n
\pm O(\ln \ln n)}{n}   ,\label{pe-lnn-lnlnn-pm}
\end{align}
then with $\kappa_v$ denoting the vertex connectivity and $\delta$
denoting the minimum vertex degree, it holds for constant integer
$\ell$ that
\begin{align}
\mathbb{P}[\hspace{2pt}(\kappa_v = \ell) \cap (\delta > \ell
)\hspace{2pt}] &= o(1).  \label{mnd_kcon}
\end{align}
\end{lemma}

We detail the proof of Lemma \ref{lem:mvd:con} below.

\subsection{The proof of Lemma \ref{lem:mvd:con}.}

For graph $G_s(n,K_n,P_n)$, let the set of vertices be $
\mathcal{V}_n = \{ v_1, v_2, \ldots , v_n \}$. Also, for
$i=1,2,\ldots, n$, we let $S_i$ denote the set of items on vertex
$v_i$. We introduce event ${\mathcal{E} (\boldsymbol{J})}$ in the
following manner:
\begin{align}
{\mathcal{E} (\boldsymbol{J})}  &= \bigcup_{\begin{subarray}{c} T
\subseteq \mathcal{V}_n, \\ |T| \geq 2.\end{subarray}} ~ \left[
\hspace{2pt} |\cup_{v_i \in T} S_i|~\leq~{J}_{ |T|}
\hspace{2pt}\right], \label{eq:E_n_defnex}
\end{align}
where $\boldsymbol{J} =[{J}_{2} , {J}_{3}, \ldots, {J}_{n } ]$ is an
$(n-1)$-dimensional integer valued array, with $J_{i}$ defined
through
\begin{align}
 J_{i} &=
\begin{cases}
\max\{ \left \lfloor (1+\varepsilon_1) K_n \right \rfloor ,
\left \lfloor \lambda_1 K_n i \right \rfloor \},~~i=2,\ldots, r_n,\\
 \left \lfloor\mu_1 P_n \right \rfloor,~~~~~~~~~~~~~~~~~~~~~~~~~i=r_n+1, \ldots, n,
\end{cases} \label{olp_xjdef}
\end{align}
for an arbitrary constant $0<\varepsilon_1<1$ and some positive
constants $\lambda_1, \mu_1$ in Lemma \ref{prp:EJ} below, where $r_n
:= \min \big( {\big \lfloor \frac{P_n}{K_n} \big \rfloor},   \big
\lfloor \frac{n}{2} \big \rfloor \big)$.

By a crude bounding argument, we get
\begin{align}
 \bP{ \hspace{2pt}(\kappa_v = {\ell}) \cap (\delta > {\ell}
)\hspace{2pt}}~~~~~~~~~~~~~~~~~~~~~~~~~~~~~ \\
\nonumber \leq \bP{\hspace{2pt}{\mathcal{E} (\boldsymbol{J})}} +
\bP{ (\kappa_v = {\ell})  \cap (\delta > {\ell} )
\cap\overline{\mathcal{E} (\boldsymbol{J})} \hspace{2pt}}.
\end{align}
Hence, a proof of Lemma \ref{lem:mvd:con} consists of proving two
lemmas below. Under (\ref{thm:uni:eq:edge}) with $|\alpha_n| = O(\ln
\ln n)$, we have $\frac{1}{s!} \cdot \frac{{K_n}^{2s}}{{P_n}^{s}} =
\frac{\ln  n \pm O(\ln \ln n)}{n} = o(1)$ and $\frac{K_n}{P_n} =
o(1)$, enabling us to use Lemmas \ref{prp:EJ} and
\ref{prp-kvl-del-EJ}.

\begin{lemma}[\hspace{-.5pt}{\cite[Proposition 3]{ZhaoYaganGligor} (our work)}\hspace{0pt}] \label{prp:EJ}
If $ P_n = \Omega(n)$ and $\frac{K_n}{P_n} = o(1)$, then for an
arbitrary constant $0<\varepsilon_1<1$ and some selected positive
constants $\lambda_1,\mu_1$, it holds that
\begin{align}
 \bP{\mathcal{E} (\boldsymbol{J})} = o(1).
\label{eq:OneLawAfterReductionPart1}
\end{align}
\end{lemma}

\begin{lemma} \label{prp-kvl-del-EJ}
For uniform random $s$-intersection graph $G_s(n,K_n,P_n)$ under
$P_n = \Omega(n)$ and
\begin{align}
\frac{1}{s!} \cdot \frac{{K_n}^{2s}}{{P_n}^{s}}   &  = \frac{\ln  n
\pm O(\ln \ln n)}{n} , \label{cn-qn-2lnn-n}
\end{align}
then
\begin{align}
  \bP{ (\kappa_v = {\ell})  \cap (\delta >
{\ell} ) \cap\overline{\mathcal{E} (\boldsymbol{J})} \hspace{2pt}} =
o(1). \label{eq:OneLawAfterReductionPart2}
\end{align}
\end{lemma}

The proof of Lemma \ref{prp-kvl-del-EJ} is given in Section
\ref{sec:prf:prop:OneLawAfterReductionPart2} in the Appendix.

\section{Establishing Theorem \ref{thm:bin}} \label{sec:prf:thm:bin}

Similar to the idea of confining $|\alpha_n |$ in Theorem
\ref{thm:uni}, here we confine $|\beta_n |$ as $O(\ln \ln n)$ in
Theorem \ref{thm:bin}. Specifically, we will demonstrate
\begin{align}
\textrm{Theorem \ref{thm:bin} under }|\beta_n | = O(\ln \ln n)
\Rightarrow  \textrm{Theorem \ref{thm:bin}}. \label{cp_alph_bin}
\end{align}
Since $k$-vertex-connectivity, $k$-edge-connectivity, and the
property of minimum vertex degree being at least $k$, are all
monotone increasing, then to show (\ref{cp_alph_bin}), it suffices
to prove the following lemma.

\begin{lemma} \label{graph_Hs_cpln}

\textbf{(a)} For graph $H_s(n,t_n,P_n)$ under
\begin{align}
\frac{1}{s!} \cdot {t_n}^{2s}{P_n}^{s}  &  = \frac{\ln  n + {(k-1)}
\ln \ln n + {\beta_n}}{n}  \label{al0-parta-Hs-od}
\end{align}
with $\lim_{n \to \infty}\beta_n = -\infty$, there exists graph
$H_s(n,\widetilde{t_n}, \widetilde{P_n})$ under
\begin{align}
\frac{1}{s!} \cdot {\widetilde{t_n}}^{2s}{\widetilde{P_n}}^{s}  &  =
\frac{\ln  n + {(k-1)} \ln \ln n + {\widetilde{\beta_n}}}{n}
\label{al0-parta-Hs}
\end{align}
with $\lim_{n \to \infty}\widetilde{\beta_n} = -\infty$ and
$\widetilde{\beta_n} = -O(\ln \ln n)$ such that there exists a graph
coupling under which $H_s(n,t_n,P_n)$ is a spanning subgraph of
$H_s(n,\widetilde{t_n},\widetilde{P_n}) $.

\textbf{(b)} For graph $H_s(n,t_n,P_n)$ under
\begin{align}
\frac{1}{s!} \cdot {t_n}^{2s}{P_n}^{s}  &  = \frac{\ln  n + {(k-1)}
\ln \ln n + {\beta_n}}{n}  \label{al0-parta-Hs-pb-od}
\end{align}
with $\lim_{n \to \infty}\beta_n = \infty$, there exists graph
$H_s(n,\widehat{t_n}, \widehat{P_n})$ under
\begin{align}
\frac{1}{s!} \cdot {\widehat{t_n}}^{2s}{\widehat{P_n}}^{s}  &  =
\frac{\ln  n + {(k-1)} \ln \ln n + {\widehat{\beta_n}}}{n}
\label{al0-parta-Hs-pb}
\end{align}
with $\lim_{n \to \infty}\widehat{\beta_n} = \infty$ and
$\widehat{\beta_n} = O(\ln \ln n)$ such that there exists a graph
coupling under which $H_s(n,t_n,P_n)$ is a spanning supergraph of
$H_s(n,\widehat{t_n},\widehat{P_n})$.

\end{lemma}

The proof of Lemma \ref{graph_Hs_cpln} is detailed in Section
\ref{sec:pro:graph_Hs_cpln} in the Appendix.

Now we use Theorem \ref{thm:uni} to prove Theorem \ref{thm:bin} with
confined $|\beta_n | = O(\ln \ln n)$. Here, the main idea is to
exploit a {\em coupling} result between the uniform $s$-intersection
graph and a binomial $s$-intersection graph. Let $\mathcal {I}^{*}$
denote either one of the following graph properties:
$k$-vertex-connectivity, $k$-edge-connectivity, and the property of
minimum vertex degree being at least $k$. With $K_n^{-}$ and
$K_n^{+}$ defined by\vspace{-5pt}
\begin{align}
K_n^{\pm} &  =
 t_n P_n \pm \sqrt{3\ln n (\ln n + t_n P_n)}
       \label{KnminandKnplus},\vspace{-5pt}
\end{align}
we have from Lemma \ref{lem:cp} that if $t_n P_n = \omega(\ln n )$,
\vspace{-5pt} then
\begin{align}
\label{propIud} & \mathbb{P} \big[\hspace{2pt}\textrm{Graph
}G_s(n,K_n^{-},P_n)\textrm{
has $\mathcal {I}^{*}$}.\hspace{2pt}\big] - o(1)  \\
&  \leq \mathbb{P} \big[\hspace{2pt}\textrm{Graph
}H_s(n,t_n,P_n)\textrm{ has $\mathcal {I}^{*}$}.
  \hspace{2pt}\big] \nonumber \\
&  \leq  \mathbb{P} \big[\hspace{2pt}\textrm{Graph
}G_s(n,K_n^{+},P_n)\textrm{ has $\mathcal
{I}^{*}$}.\hspace{2pt}\big] \vspace{-5pt} + o(1).\nonumber
 \end{align}

 Under conditions (\ref{thm:bin:eq:P}), (\ref{thm:bin:eq:edge}),
 and $|\beta_n | = O(\ln \ln n)$, we now show that
 $t_n P_n = \omega(\ln n )$. From (\ref{thm:bin:eq:edge})
 and $|\beta_n | = O(\ln \ln n)$, we first get\vspace{-5pt}
\begin{align}
\frac{1}{s!} \cdot {t_n}^{2s}{P_n}^{s}  & =   \frac{\ln  n \pm O(\ln
\ln n) }{n} = \frac{\ln n}{n}  \cdot [1\pm o(1)] \vspace{-5pt} .
\label{tn2sPnslnneq}
\end{align}
From (\ref{thm:bin:eq:P}) and (\ref{tn2sPnslnneq}), it follows
that\vspace{-5pt}
\begin{align}
&\label{eq_tnPn_c} t_n P_n  = \sqrt{ {t_n}^{2}{P_n}} \cdot \sqrt{P_n
}  \\ & =
\begin{cases}
 \big\{ {s! n^{-1} \ln n} \cdot [1\pm o(1)] \big\}^{\frac{1}{2s}} \cdot \sqrt{\Omega(n)}  ,   &\textrm{for } s \geq 2, \\
 \big\{  {s! n^{-1} \ln n}  \cdot [1\pm o(1)]  \big\}^{\frac{1}{2s}} \cdot \sqrt{\Omega(n^c)}  ,   &\textrm{for } s = 1,
 \end{cases} \nonumber \\ &
=
\begin{cases}
  \Omega\big(n^{\frac{1}{2}-\frac{1}{2s}} (\ln n)^{\frac{1}{2s}}\big) , &\textrm{for } s \geq 2, \\
  \Omega\big(n^{\frac{c-1}{2}} (\ln n)^{\frac{1}{2}} \big), &\textrm{for } s = 1,\vspace{-5pt}
 \end{cases}\nonumber
\end{align}
yielding $t_n P_n = \omega(\ln n )$ in view of $c>1$, so we can use
(\ref{propIud}).
 Using (\ref{KnminandKnplus}) and (\ref{eq_tnPn_c}),
we further\vspace{-5pt} obtain
\begin{align}
\label{KnsPntn2s} \frac{{(K_n^{\pm})}^{2s}}{{P_n}^{s}}  &  =
  \frac{{\big[ t_n P_n \pm \sqrt{3\ln n (\ln n + t_n P_n)}\hspace{1pt}\big]}^{2s}}{{P_n}^{s}}
\\ & = \frac{(t_n P_n )^{2s}}{{P_n}^{s}} \cdot \Bigg[1 \pm  \sqrt{\frac{3\ln n}{t_n P_n} \bigg(\frac{\ln n}{t_n P_n} + 1\bigg)} \hspace{2pt} \Bigg]^{2s}
    \nonumber \\ & = {t_n}^{2s} {P_n}^{s}  \cdot \bigg[1 \pm o\bigg(\frac{1}{\ln n}\bigg)\bigg], \nonumber
\end{align}
where in the last step we use $t_n P_n = \omega\big((\ln n
)^3\big)$, which follows from (\ref{eq_tnPn_c}) due to constant
$c>1$.

  Applying (\ref{thm:bin:eq:edge}) and (\ref{tn2sPnslnneq}) to (\ref{KnsPntn2s}), we have
 \begin{align}
\frac{1}{s!}  \cdot  \frac{{(K_n^{\pm })}^{2s}}{{P_n}^{s}}  &  =
  \frac{\ln  n + {(k-1)} \ln \ln n + {\beta_n \pm o(1)}}{n}. \label{KnsPntn2sbtn}
\end{align}
In view of (\ref{KnsPntn2sbtn}) and $P_n = \Omega(n )$, we use
Theorem \ref{thm:uni} to obtain
\begin{align}
\label{KnsPntn2sbtnud} &\lim_{n \to \infty} \mathbb{P}
\big[\hspace{2pt}\textrm{Graph }G_s(n,K_n^{\pm },P_n)\textrm{ has
$\mathcal {I}^{*}$}.\hspace{2pt}\big]  \\ & = e^{- \frac{e^{-\lim_{n
\to \infty} [\beta_n \pm o(1)]}}{(k-1)!}} = e^{- \frac{e^{-\lim_{n
\to \infty}\beta_n}}{(k-1)!}}.\nonumber
 \end{align}

 The proof of Theorem \ref{thm:bin} is completed
 by (\ref{propIud}) and (\ref{KnsPntn2sbtnud}).

 \section{Conclusion} \label{sec:Conclusion}

Random $s$-intersection graphs have been used in a wide range of
applications. Two extensively studied models are a uniform random
$s$-intersection graph and a binomial random $s$-intersection graph.
In this paper, for a uniform/binomial random $s$-intersection graph,
we derive exact asymptotic expressions for the probabilities of the
following three properties: (i) $k$-vertex-connectivity, (ii)
$k$-edge-connectivity and (iii) the property that each vertex has
degree at least $k$.


\section*{Acknowledgements}

This research was supported in part by CMU CyLab under grant
CNS-0831440 from the National Science Foundation. The views and
conclusions contained in this document are those of the authors and
should not be interpreted as representing the official policies,
either expressed or implied, of any sponsoring institution, the U.S.
government or any other entity.


{


\section{Appendix}

We first present in Section \ref{sec:add:lemma} additional lemmas
used in proving the theorems. Afterwards, we detail the proofs of
the lemmas.

\subsection{Additional lemmas.} \label{sec:add:lemma}

Some additional lemmas are given below. The relation ``$\sim$''
stands for an asymptotical equivalence; i.e., $f_n \sim g_n$ means
$\lim_{n \to
  \infty}({f_n }/{g_n })=1$.

\begin{lemma}  \label{lem:logn2}
If $\frac{1}{s!} \cdot \frac{{K_n}^{2s}}{{P_n}^{s}}    = \frac{\ln
n \pm O(\ln \ln n) }{n}$ and $P_n = \Omega(n^c)$ for constant $c$,
then $K_n = \Omega\big( n^{\frac{c}{2} - \frac{1}{2s}} (\ln
n)^{\frac{1}{2s}} \big)$.
\end{lemma}

\begin{lemma} \label{qn-dist}

The following properties (a) and (b) hold, where $q_n$ is the edge
probability in uniform random $s$-intersection graph
$G_s(n,K_n,P_n)$.

\begin{itemize}[leftmargin=20pt]
\item[(a)] If $P_n = \Omega(n)$ and $\frac{1}{s!} \cdot \frac{{K_n}^{2s}}{{P_n}^{s}}   = \frac{\ln  n \pm O(\ln \ln n) }{n}$, then $q_n  \sim \frac{1}{s!} \cdot \frac{{K_n}^{2s}}{{P_n}^{s}}$ and $\big| \hspace{2pt} q_n - \frac{1}{s!} \cdot \frac{{K_n}^{2s}}{{P_n}^{s}} \hspace{2pt} \big|     =   o\big(\frac{1}{n}\big)$.
\item[(b)] If $P_n = \Omega(n)$ and $q_n = \frac{\ln  n \pm O(\ln \ln n) }{n}$, then $q_n  \sim \frac{1}{s!} \cdot \frac{{K_n}^{2s}}{{P_n}^{s}}$ and $\big| \hspace{2pt} q_n - \frac{1}{s!} \cdot \frac{{K_n}^{2s}}{{P_n}^{s}} \hspace{2pt} \big|     =   o\big(\frac{1}{n}\big)$.
\end{itemize}

\end{lemma}

\begin{lemma} \label{lem_prob_Eij_S1r}

For uniform random $s$-intersection graph $G_s(n,K_n,P_n)$ under
$K_n = \omega(1)$, the following properties (a) (b) and (c) hold for
$i = r+1, r+2, \ldots, n$ (i.e., vertex $v_i \notin \{v_1, v_2,
\ldots, v_r\}$), where $E_{ij}$ denotes the event that an edge
exists between vertices $v_i$ and $v_j$, $S_i$ is the number of
items on vertex $v_i$, and $q_n$ is the edge probability.
\begin{itemize}[leftmargin=20pt]
\item[(a)] If $|\bigcup_{j=1}^{r} S_j| \geq \lfloor  (1+{\varepsilon_1}) K_n \rfloor$\vspace{2pt} for some positive constant $\varepsilon_1$, then
for any positive constant $\varepsilon_2 < (1+{\varepsilon_1})^s -
1$, it holds for all $n$ sufficiently large that
\begin{align}
 \mathbb{P}\bigg[\hspace{2pt}\bigcap_{j=1}^{r}\overline{E_{ij}} \hspace{2pt} \bigg| \hspace{2pt} S_1, S_2, \ldots, S_r\bigg]
& \leq  e^{- q_n (1+\varepsilon_2)}.
\end{align}
\item[(b)] If $|\bigcup_{j=1}^{r} S_j| \geq \lfloor  \lambda_1 r K_n \rfloor$\vspace{2pt} for some positive constant $\lambda_1$, then
for any positive constant $\lambda_2 < {\lambda_1}^s$, it holds for
all $n$ sufficiently large that
\begin{align}
 \mathbb{P}\bigg[\hspace{2pt}\bigcap_{j=1}^{r}\overline{E_{ij}} \hspace{2pt} \bigg| \hspace{2pt} S_1, S_2, \ldots, S_r\bigg]
& \leq  e^{- \lambda_2 r q_n}.
\end{align}
\item[(c)] If $|\bigcup_{j=1}^{r} S_j| \geq \lfloor \mu_1 P_n \rfloor$\vspace{2pt} for some positive constant $\mu_1$,
then for any positive constant $\mu_2 < (s!)^{-1}{\mu_1}^s$, it
holds for all $n$ sufficiently large that
\begin{align}
 \mathbb{P}\bigg[\hspace{2pt}\bigcap_{j=1}^{r}\overline{E_{ij}} \hspace{2pt} \bigg| \hspace{2pt} S_1, S_2, \ldots, S_r\bigg]
& \leq  e^{- \mu_2 K_n}.
\end{align}
\end{itemize}

\end{lemma}

\begin{lemma} \label{olp_lem1}

For uniform random $s$-intersection graph $G_s(n,K_n,P_n)$ under
$P_n = \Omega(n)$, $K_n = \omega(1)$ and $r_n := \min \big( {\big
\lfloor \frac{P_n}{K_n} \big \rfloor}, \big \lfloor \frac{n}{2} \big
\rfloor \big) = \omega(1)$, the following properties (a) (b) and (c)
hold for any constant integer $R \geq 2$, where $\varepsilon_1$,
$\lambda_1$ and $\mu_1$ are specified in Lemma \ref{prp:EJ}, and events $\mathcal{A}_{{\ell},r}$ and $\mathcal{E} (\boldsymbol{J})$ are defined in Sections 6.3 and 3.3, respectively.

\begin{itemize}[leftmargin=20pt]
\item[(a)] Let $ \varepsilon_3$ be any positive constant with $ \varepsilon_3 < (1+\varepsilon_1)^s-1$. For all $n$ sufficiently large, it holds for $r = 2, 3, \ldots, R$ that
\begin{align}
 \bP{ \mathcal{A}_{{\ell},r} \hspace{1.5pt}\cap\hspace{1.5pt} \overline{\mathcal{E} (\boldsymbol{J})} }  &\leq
 r^{r-2} {q_n}^{r-1}   ( r {q_n})^{{\ell} }   e^{- q_n n (1+\varepsilon_3)}
 . \nonumber
\end{align}
\item[(b)] Let $\lambda_2$ be any positive constant with $\lambda_2 < {\lambda_1}^s$.
 For all $n$ sufficiently large, it holds for $r = R + 1, R + 2, \ldots,  r_n $ that
\begin{align}
 \bP{ \mathcal{A}_{{\ell},r} \hspace{2pt}\cap\hspace{2pt} \overline{\mathcal{E} (\boldsymbol{J})} }  &\leq
 r^{r-2} {q_n}^{r-1}   e^{- \lambda_2 r q_n n /3}
 . \nonumber
\end{align}
\item[(c)] Let $\mu_2$ be any positive constant with $\mu_2 < (s!)^{-1}{\mu_1}^s$.
 For all $n$ sufficiently large, it holds for $r = r_n + 1, r_n + 2, \ldots, \lfloor
\frac{n-{\ell}}{2} \rfloor $ that
\begin{align}
 \bP{ \mathcal{A}_{{\ell},r} \hspace{2pt}\cap\hspace{2pt} \overline{\mathcal{E} (\boldsymbol{J})} }  &\leq
 e^{- \mu_2  K_n  n /3}
 . \nonumber
\end{align}
\end{itemize}

\end{lemma}

\begin{lemma}[\hspace{-.1pt}{\cite[Lemma 4]{Rybarczyk}}]  \label{lem:cp}

Let $K_n^{-}$ and $K_n^{+}$ denote $t_n P_n - \sqrt{3\ln n (\ln n +
t_n P_n)}$ and $t_n P_n + \sqrt{3\ln n (\ln n + t_n P_n)}$,
respectively. If $t_n P_n = \omega(\ln n )$, then for any monotone
increasing graph property $\mathcal {I}$, it holds that
\begin{align}
 & \mathbb{P} \big[\hspace{2pt}\textrm{Graph }G_s(n,K_n^{-},P_n)\textrm{
has $\mathcal {I}$}.\hspace{2pt}\big] - o(1) \nonumber \\
&  \leq \mathbb{P} \big[\hspace{2pt}\textrm{Graph
}H_s(n,t_n,P_n)\textrm{ has $\mathcal {I}$}.
  \hspace{2pt}\big] \nonumber \\
&  \leq  \mathbb{P} \big[\hspace{2pt}\textrm{Graph
}G_s(n,K_n^{+},P_n)\textrm{ has $\mathcal {I}$}.\hspace{2pt}\big]  +
o(1). \nonumber
 \end{align}

\end{lemma}

\subsection{Proof of Lemma \ref{graph_Gs_cpl}} \label{sec_graph_Gs_cpl}

\paragraph{Proving property (a).}

 We
 define $\widetilde{\alpha_n}^*$ by
 \begin{align}
\widetilde{\alpha_n}^* &  = \max\{\alpha_n, -\ln \ln n\},
\label{al2-parta}
\end{align}
and define $\widetilde{K_n}^*$ such that
\begin{align}
\frac{1}{s!} \cdot \frac{({\widetilde{K_n}^{*}})^{2s}}{{{P_n}}^{s}}
&  = \frac{\ln  n + {(k-1)} \ln \ln n + \widetilde{\alpha_n}^*}{n}.
\label{al3-parta}
\end{align}
We set
\begin{align}
\widetilde{K_n} & : = \big\lfloor \widetilde{K_n}^* \big\rfloor,
\label{al4-parta}
\end{align}
and
\begin{align}
\widetilde{P_n} & : = P_n.  \label{al5-parta}
\end{align}

From (\ref{al1}) (\ref{al2-parta}) and (\ref{al3-parta}), it holds
that
\begin{align}
K_n \leq \widetilde{K_n}^*.  \label{Kn1-parta}
\end{align}
Then by (\ref{al4-parta}) (\ref{Kn1-parta}) and the fact that $K_n$
and $\widetilde{K_n}$ are both integers, it follows that
\begin{align}
K_n \leq \widetilde{K_n}.  \label{al6-parta}
\end{align}
 From (\ref{al5-parta}) and (\ref{al6-parta}), by  \cite[Lemma 3]{Rybarczyk}, there exists a graph coupling under which $G_s(n,K_n,P_n)$ is a spanning subgraph of $G_s(n,\widetilde{K_n},\widetilde{P_n})$. Therefore, the proof of property (a) is completed once we show
$\widetilde{\alpha_n}$ defined in $(\ref{al0-parta})$ satisfies
\begin{align}
  \lim_{n \to \infty}\widetilde{\alpha_n} & = - \infty, \label{al8-parta} \\
 \widetilde{\alpha_n} & = - O(\ln \ln n).  \label{al7-parta}
\end{align}

We first prove (\ref{al8-parta}). From (\ref{al0-parta})
(\ref{al3-parta}) and (\ref{al4-parta}), it holds that
\begin{align}
\widetilde{\alpha_n} \leq \widetilde{\alpha_n}^*, \label{haa-parta}
\end{align}
which together with (\ref{al2-parta}) and $\lim_{n \to
\infty}\alpha_n = -\infty$ yields (\ref{al8-parta}).

Now we establish (\ref{al7-parta}). From (\ref{al4-parta}), we have
$\widetilde{K_n} > \widetilde{K_n}^* - 1$. Then from
(\ref{al0-parta}) and (\ref{al5-parta}), it holds that
\begin{align}
 \label{aph1-parta}\widetilde{\alpha_n} = n \cdot \frac{1}{s!} \cdot \frac{{\widetilde{K_n}}^{2s}}{{{P_n}}^{s}}  - [\ln  n + {(k-1)} \ln \ln n]~~~  \\
   >  n \cdot \frac{1}{s!} \cdot \frac{{(\widetilde{K_n}^* - 1)}^{2s}}{{{P_n}}^{s}}  - [\ln  n + {(k-1)} \ln \ln n]  .\nonumber
\end{align}
By $\lim_{n \to \infty}\alpha_n =- \infty$, it holds that $\alpha_n
\leq 0$ for all $n$ sufficiently large. Then from (\ref{al2-parta}),
it follows that
\begin{align}
 \widetilde{\alpha_n}^* = - O(\ln \ln n),  \label{widetilde-al2-parta}
\end{align}
which along with Lemma \ref{lem:logn2}, equation (\ref{al3-parta})
and condition $P_n = \Omega(n)$ induces
\begin{align}
 \widetilde{K_n}^* & = \Omega\big( (\ln n)^{\frac{1}{2s}} \big). \label{aph5-parta}
\end{align}
Hence, we have $\lim_{n \to \infty} \widetilde{K_n}^* = \infty$ and
it further holds for all $n$ sufficient large that
\begin{align}
{(\widetilde{K_n}^* - 1)}^{2s} > ({\widetilde{K_n}^{*}})^{2s} - 3s
({\widetilde{K_n}^{*}})^{2s-1}. \label{aph2-parta}
\end{align}
Applying (\ref{aph2-parta}) to (\ref{aph1-parta}) and then using
(\ref{al3-parta}), Lemma \ref{lem:logn2} and $P_n = \Omega(n)$, it
follows that
\begin{align}
 \label{widetilde-al-parta}  &   \widetilde{\alpha_n} \\
  &   >  \frac{n}{s!}\hspace{-1pt} \cdot\hspace{-1pt}
   \frac{ ({\widetilde{K_n}^{*}})^{2s}
 \hspace{-1pt}-\hspace{-1pt} 3s  ({\widetilde{K_n}^{*}})^{2s-1}}{{{P_n}}^{s}}
  \hspace{-1pt}-\hspace{-1pt} [\ln  n \hspace{-1pt}+\hspace{-1pt} {(k\hspace{-1pt}-\hspace{-1pt}1)} \ln \ln n]  \nonumber \\
 &    = \widetilde{\alpha_n}^* - \frac{3s}{s!} \cdot n \cdot \Theta\big({P_n}^{-\frac{1}{2}} n^{-\frac{2s-1}{2s}} (\ln n)^{\frac{2s-1}{2s}} \big)   \nonumber \\
 &    = \widetilde{\alpha_n}^* - O\big(n^{-\frac{1}{2}+\frac{1}{2s}} (\ln n)^{1-\frac{1}{2s}}\big).\nonumber
\end{align}

As noted at the beginning of Section \ref {sec:basic:ideas}, our
proof is for $s \geq 2$ since the case of $s = 1$ already is proved
by us \cite{ZhaoCDC}. Using $s \geq 2$ in
(\ref{widetilde-al-parta}), it holds that $ \widetilde{\alpha_n} >
\widetilde{\alpha_n}^* + o(1)$, which along with (\ref{haa-parta})
and (\ref{widetilde-al2-parta}) yields (\ref{al7-parta}).

\paragraph{Proving property (b).}

 We
 define $\widehat{\alpha_n}^*$ by
 \begin{align}
\widehat{\alpha_n}^* &  = \min\{\alpha_n, \ln \ln n\}, \label{al2}
\end{align}
and define $\widehat{K_n}^*$ such that
\begin{align}
\frac{1}{s!} \cdot \frac{({\widehat{K_n}^{*}})^{2s}}{{{P_n}}^{s}}  &
= \frac{\ln  n + {(k-1)} \ln \ln n + \widehat{\alpha_n}^*}{n}.
\label{al3}
\end{align}
We set
\begin{align}
\widehat{K_n} & : = \big\lceil \widehat{K_n}^* \big\rceil,
\label{al4}
\end{align}
and
\begin{align}
\widehat{P_n} & : = P_n.  \label{al5}
\end{align}

From (\ref{al1}) (\ref{al2}) and (\ref{al3}), it holds that
\begin{align}
K_n \geq \widehat{K_n}^*.  \label{Kn1}
\end{align}
Then by (\ref{al4}) (\ref{Kn1}) and the fact that $K_n$ and
$\widehat{K_n}$ are both integers, it follows that
\begin{align}
K_n \geq \widehat{K_n}.  \label{al6}
\end{align}
 From (\ref{al5}) and (\ref{al6}), by  \cite[Lemma 3]{Rybarczyk}, there exists a graph coupling under which $G_s(n,K_n,P_n)$ is a spanning supergraph of $G_s(n,\widehat{K_n},\widehat{P_n})$. Therefore, the proof of property (b) is completed once we show
$\widehat{\alpha_n}$ defined in $(\ref{al0})$ satisfies
\begin{align}
  \lim_{n \to \infty}\widehat{\alpha_n} & = \infty, \label{al8} \\
 \widehat{\alpha_n} & = O(\ln \ln n).  \label{al7}
\end{align}

We first prove (\ref{al8}). From (\ref{al0}) (\ref{al3}) and
(\ref{al4}), it holds that
\begin{align}
\widehat{\alpha_n} \geq \widehat{\alpha_n}^*, \label{haa}
\end{align}
which together with (\ref{al2}) and $\lim_{n \to \infty}\alpha_n =
\infty$ yields (\ref{al8}).

Now we establish (\ref{al7}). From (\ref{al4}), we have
$\widehat{K_n} < \widehat{K_n}^* + 1$. Then from (\ref{al0}) and
(\ref{al5}), it holds that
\begin{align}
\label{aph1} \widehat{\alpha_n}  = n \cdot \frac{1}{s!} \cdot
\frac{{\widehat{K_n}}^{2s}}
{{{P_n}}^{s}}  - [\ln  n + {(k-1)} \ln \ln n]~~~  \\
 <  n \cdot \frac{1}{s!}
  \cdot \frac{{(\widehat{K_n}^* + 1)}^{2s}}{{{P_n}}^{s}}
   - [\ln  n + {(k-1)} \ln \ln n]  .\nonumber
\end{align}
By $\lim_{n \to \infty}\alpha_n = \infty$, it holds that $\alpha_n
\geq 0$ for all $n$ sufficiently large. Then from (\ref{al2}), it
follows that
\begin{align}
 \widehat{\alpha_n}^* = O(\ln \ln n),  \label{widehat-al2}
\end{align}
which along with Lemma \ref{lem:logn2}, equation (\ref{al3}) and
condition $P_n = \Omega(n)$ induces
\begin{align}
 \widehat{K_n}^* & = \Omega\big( (\ln n)^{\frac{1}{2s}} \big). \label{aph5}
\end{align}
Hence, we have $\lim_{n \to \infty} \widehat{K_n}^* = \infty$ and it
further holds for all $n$ sufficient large that
\begin{align}
{(\widehat{K_n}^* + 1)}^{2s}< ({\widehat{K_n}^{*}})^{2s} + 3s
({\widehat{K_n}^{*}})^{2s-1}. \label{aph2}
\end{align}
Applying (\ref{aph2}) to (\ref{aph1}) and then using (\ref{al3}),
Lemma \ref{lem:logn2} and $P_n = \Omega(n)$, it follows that
\begin{align}
\label{widehat-al} & \widehat{\alpha_n} \\
 \nonumber  & <  \frac{n}{s!} \hspace{-1pt}\cdot\hspace{-1pt} \frac{ ({\widehat{K_n}^{*}})^{2s}\hspace{-1pt} +\hspace{-1pt} 3s  ({\widehat{K_n}^{*}})^{2s-1}}{{{P_n}}^{s}}  \hspace{-1pt}-\hspace{-1pt} [\ln  n \hspace{-1pt}+\hspace{-1pt} {(k\hspace{-1pt}-\hspace{-1pt}1)} \ln \ln n]  \nonumber \\
 &    = \widehat{\alpha_n}^* + \frac{3s}{s!} \cdot n \cdot \Theta\big({P_n}^{-\frac{1}{2}} n^{-\frac{2s-1}{2s}} (\ln n)^{\frac{2s-1}{2s}} \big)   \nonumber \\
 &    = \widehat{\alpha_n}^* + O\big(n^{-\frac{1}{2}+\frac{1}{2s}} (\ln
 n)^{1-\frac{1}{2s}}\big). \nonumber
\end{align}

As noted at the beginning of Section \ref {sec:basic:ideas}, our
proof is for $s \geq 2$ since the case of $s = 1$ already is proved
by Rybarczyk \cite{ryb3}. Using $s \geq 2$ in (\ref{widehat-al}), it
holds that $ \widehat{\alpha_n} < \widehat{\alpha_n}^* + o(1)$,
which along with (\ref{haa}) and (\ref{widehat-al2}) yields
(\ref{al7}).

 \subsection{The proof of Lemma \ref{prp-kvl-del-EJ}.} \label{sec:prf:prop:OneLawAfterReductionPart2}

By the analysis in \cite[Section IV]{ZhaoYaganGligor}, we obtain
\cite[Equation (148)]{ZhaoYaganGligor}. Namely, with some events
defined as follows:
\begin{itemize}[leftmargin=20pt]
\item $\mathcal{C}_{r}$: event that the induced subgraph of $G_s(n, K_n, P_n)$ defined on vertex set $\{v_1, v_2, \ldots, v_r\}$ is connected,
\item $\mathcal{B}_{\ell,r}$: event that any vertex in $\{v_{r+1}, v_{r+2}, \ldots, v_{r+\ell} \}$ has an edge with at least one vertex in $\{v_{1}, v_{2}, \ldots,  v_{r}\}$,
\item $\mathcal{D}_{\ell,r}$: event that any vertex in $\{v_{r+\ell+1}, v_{r+\ell+2}, \ldots, v_{n} \}$ and any vertex in $\{v_{1}, v_{2}, \ldots,  v_{r}\}$ has no edge in between, and
\item  $\mathcal{A}_{\ell, r}$: event that events $\mathcal{C}_{r}$, $\mathcal{B}_{\ell,r}$ and $\mathcal{D}_{\ell,r}$ all happen,
\end{itemize}
it holds that
\begin{align}
 \label{eq:BasicIdea+UnionBound2} \bP{(\kappa ={\ell}) ~\cap~ (\delta > {\ell}) ~\cap~
\overline{\mathcal{E} (\boldsymbol{J})} }~~~~~~~~~~~~~  \\
 \leq   \sum_{r=2}^{ \lfloor \frac{n-{\ell}}{2} \rfloor } {n \choose
{\ell} }{ {n-{\ell}} \choose r} ~ \bP{ \mathcal{A}_{{\ell},r} ~\cap~
\overline{\mathcal{E} (\boldsymbol{J})}} . \nonumber
\end{align}

The proof of Lemma \ref{prp-kvl-del-EJ} is completed once we show
the following three results:
\begin{align}
\sum_{r=2}^{ R} {n \choose {\ell} }{ {n-{\ell}} \choose r} ~ \bP{
\mathcal{A}_{{\ell},r} ~\cap~ \overline{\mathcal{E}
(\boldsymbol{J})}} = o(1) , \label{prf-e1}
\\  \sum_{r=R+1}^{ r_n} {n \choose {\ell} }{ {n-{\ell}} \choose
r} ~ \bP{ \mathcal{A}_{{\ell},r} ~\cap~ \overline{\mathcal{E}
(\boldsymbol{J})}} = o(1)  , \label{prf-e2}
\end{align}
and
\begin{align}
\sum_{r=r_n+1}^{  \lfloor \frac{n-{\ell}}{2}\rfloor } {n \choose
{\ell} }{ {n-{\ell}} \choose r} ~ \bP{ \mathcal{A}_{{\ell},r} ~\cap~
\overline{\mathcal{E} (\boldsymbol{J})}} = o(1)  , \label{prf-e3}
\end{align}
where $r_n = \min \big( {\big \lfloor \frac{P_n}{K_n} \big \rfloor},
\big \lfloor \frac{n}{2} \big \rfloor \big)$.

From condition (\ref{cn-qn-2lnn-n}), it follows that $
\frac{K_n}{P_n} = o(1)$, yielding $r_n = \omega(1)$. From conditions
(\ref{cn-qn-2lnn-n}) and $P_n = \Omega(n)$, we use Lemma
\ref{lem:logn2} to derive $K_n = \omega(1)$. Therefore, we have $P_n
= \Omega(n)$, $K_n = \omega(1)$ and $r_n = \omega(1)$, enabling us
to use Lemma \ref{olp_lem1}.

In addition, given conditions (\ref{cn-qn-2lnn-n}) and $P_n =
\Omega(n)$, we use Lemma \ref{qn-dist} to obtain
\begin{align}
q_n & =  \frac{\ln  n \pm O(\ln \ln n)}{n}
.\label{cn-qn-2lnn-n2-newqn}
\end{align}
Hence, it holds that
\begin{align}
q_n & \leq  \frac{2\ln  n}{n}, \textrm{ for all $n$ sufficiently
large},\label{cn-qn-2lnn-n2-newqn-1}
\end{align}
and there exists constant $c_0$ such that
\begin{align}
q_n & \geq  \frac{\ln  n - c_0 \ln \ln n}{n}, \textrm{ for all $n$
sufficiently large}.\label{cn-qn-2lnn-n2-newqn-2}
\end{align}

\subsubsection{Establishing (\ref{prf-e1}).}

From ${n \choose {\ell}} \leq n^{{\ell} }$, ${ n-{\ell} \choose r}
\leq n^{ r}$ and property (a) of Lemma \ref{olp_lem1}, it follows
that
\begin{align}
\label{olp_zja} & {n \choose {\ell} }{ {n-{\ell}} \choose r} \bP{
\mathcal{A}_{{\ell},r} ~\cap~ \overline{\mathcal{E}
(\boldsymbol{J})}}  \\ &\leq n^{{\ell} } \cdot n^{ r} \cdot r^{r-2}
{q_n}^{r-1}  ( r q_n)^{\ell} \cdot e^{-q_n n (1+\varepsilon_3) }
\nonumber \\ & =  r ^{\ell} r^{r-2} \cdot n^{{\ell}+r} {q_n}^{
{\ell}+r -1} \cdot e^{-q_n n (1+\varepsilon_3) } . \nonumber
\end{align}

Applying (\ref{cn-qn-2lnn-n2-newqn}) and
(\ref{cn-qn-2lnn-n2-newqn-2}) to (\ref{olp_zja}), we get
 \begin{align}
& {n \choose {\ell} }{ {n-{\ell}} \choose
r} \bP{ \mathcal{A}_{{\ell},r} ~\cap~ \overline{\mathcal{E} (\boldsymbol{J})}}  \nonumber \\
&\leq  r ^{\ell} r^{r-2}   n^{{\ell}+r}   \bigg(\hspace{-2pt} \frac{2 \ln  n  }{n}\hspace{-1pt}\bigg)^{ {\ell}+r -1}  \hspace{-1pt} e^{-(1+\varepsilon_3)(\ln  n - c_0 \ln \ln n)}  \nonumber \\
& \leq 2^{ {\ell}+r -1}     r ^{\ell+r-2}   n^{-\varepsilon_3}  (\ln n)^{\ell + r -1 +c_0 (1+\varepsilon_3)}  \nonumber \\
& = o(1). \nonumber
\end{align}
Since $R$ is a constant, (\ref{prf-e1}) clearly follows.

\subsubsection{Establishing (\ref{prf-e2}).}

From ${n \choose {\ell}} \leq n^{{\ell} }$, ${ n-{\ell} \choose r}
\leq \left( \frac{e (n-{\ell}) }{r}\right)^r \leq \left( \frac{e n
}{r}\right)^r$ and property (b) of Lemma \ref{olp_lem1}, we have
\begin{align}
 &\label{olpewrwfr2tae} {n \choose {\ell} }{ {n-{\ell}} \choose
r} \bP{ \mathcal{A}_{{\ell},r} ~\cap~ \overline{\mathcal{E}
(\boldsymbol{J})}} \\ &\leq  n^{\ell} \cdot \left( \frac{e
(n-{\ell}) }{r}\right)^r \cdot r^{r-2}{q_n}^{r-1}  e^{- \lambda_2 r
{q_n} n / 3}  \nonumber \\ & \leq n^{{\ell}+r} e^{r} {q_n}^{r-1}
e^{- \lambda_2 r {q_n} n / 3}.\nonumber
\end{align}

 Applying (\ref{cn-qn-2lnn-n2-newqn}) and (\ref{cn-qn-2lnn-n2-newqn-2}) to (\ref{olpewrwfr2tae}), we get
\begin{align}
 \label{olp_an3} &  {n \choose {\ell} }{ {n-{\ell}} \choose
r} \bP{ \mathcal{A}_{{\ell},r} ~\cap~ \overline{\mathcal{E}
(\boldsymbol{J})}}  \\ &\leq  n^{{\ell}+r} e^{r}
 \cdot \left(\frac{2\ln n}{n}\right)^{r-1} \cdot
   e^{- \lambda_2 r (\ln  n - c_0 \ln \ln n) / 3} \nonumber \\
 & \leq
n^{{\ell}+1} \cdot \big(2en^{-\lambda_2 /3} (\ln n)^{c_0\lambda_2/3
+ 1}\big)^r. \nonumber
\end{align}
Given $2en^{-\lambda_2 /3} (\ln n)^{c_0\lambda_2/3 + 1} = o(1)$ and
(\ref{olp_an3}), we obtain
\begin{align}
\nonumber & \sum_{r=R+1}^{  r_n } {n \choose {\ell} }{
{n-{\ell}} \choose r}  \bP{ \mathcal{A}_{{\ell},r} ~\cap~
\overline{\mathcal{E} (\boldsymbol{J})}} \\ & \leq
 \sum_{r=R+1}^{  \infty } n^{{\ell}+1} \cdot \big(2en^{-\lambda_2 /3} (\ln n)^{c_0\lambda_2/3 + 1}\big)^r
  \nonumber \\   &  = n^{{\ell}+1} \cdot
\frac{\big(2en^{-\lambda_2 /3} (\ln n)^{c_0\lambda_2/3 +
1}\big)^{R+1}}{1- 2en^{-\lambda_2 /3} (\ln n)^{c_0\lambda_2/3 + 1}}
 \nonumber \\
  &  \sim n^{{\ell}+1 -\lambda_2 (R+1) /3 }
  \big(2e (\ln n)^{c_0\lambda_2/3 + 1}\big)^{R+1} .\label{olp_an3ar}
\end{align}
 We pick constant $R \geq \frac{3({\ell}+1)}{\lambda_2}$ so that ${\ell}+1 -\lambda_2
(R+1) /3 \leq -\frac{\lambda_2}{3}$. As a result, we obtain
\begin{align}
 \textrm{R.H.S. of
(\ref{olp_an3ar})} &  =  o(1) \nonumber
\end{align}
 and thus establish (\ref{prf-e2}).
 
%

\subsubsection{Establishing (\ref{prf-e3}).}

 From ${n \choose
{\ell}} \leq n^{{\ell} }$ and property (c) of Lemma \ref{olp_lem1},
it holds that
\begin{align}
\label{sumrrrn1}&\sum_{r=r_n+1}^{  \lfloor \frac{n-{\ell}}{2}\rfloor
} {n \choose {\ell} }{ {n-{\ell}} \choose r} ~ \bP{
\mathcal{A}_{{\ell},r} ~\cap~ \overline{\mathcal{E}
(\boldsymbol{J})}} \\ &   \leq n^{\ell}  \cdot e^{- \mu_2  K_n  n
/3} \cdot \sum_{r=r_n+1}^{  \lfloor \frac{n-{\ell}}{2}\rfloor } {
{n-{\ell}} \choose r} . \nonumber
\end{align}
Given conditions $P_n = \Omega(n)$ and (\ref{cn-qn-2lnn-n}), we use
Lemma \ref{lem:logn2} to derive
\begin{align}
K_n = \Omega\big( n^{\frac{1}{2} - \frac{1}{2s}} (\ln
n)^{\frac{1}{2s}} \big) = \omega(1), \nonumber
\end{align}
which yields
\begin{align}
\mu _2 K_n  / 3 \geq 2 \ln 2, \textrm{ for all $n$ sufficiently
large}. \label{sumrrrn1c}
\end{align}
We have
\begin{align}
\sum_{r=r_n+1}^{  \lfloor \frac{n-{\ell}}{2}\rfloor } { {n-{\ell}}
\choose r}  \leq \sum_{r=r_n+1}^{  \lfloor \frac{n-{\ell}}{2}\rfloor
} { n \choose r } \leq \sum_{r=0}^{n} { n \choose r }  = 2^n.
\label{sumrrrn1b}
\end{align}

Applying (\ref{sumrrrn1c}) and  (\ref{sumrrrn1b}) to
(\ref{sumrrrn1}), we finally obtain
\begin{align}
&\sum_{r=r_n+1}^{  \lfloor \frac{n-{\ell}}{2}\rfloor } {n \choose
{\ell} }{ {n-{\ell}} \choose r} ~ \bP{ \mathcal{A}_{{\ell},r} ~\cap~
\overline{\mathcal{E} (\boldsymbol{J})}}\nonumber \\ &   \leq
n^{\ell} \cdot 2^n \cdot    e^{- \mu_2  K_n  n /3} \nonumber \\ & =
e^{\ell \ln n + n \ln 2 - \mu_2  K_n  n /3}
 \nonumber \\
 & \leq  e^{\ell \ln n - n \ln 2 }, \textrm{ for all $n$ sufficiently large}.\nonumber
\end{align}
The result (\ref{prf-e3}) clearly follows with $n \to \infty$.

\subsection{Proof of Lemma \ref{graph_Hs_cpln}} \label{sec:pro:graph_Hs_cpln}~

\noindent \textbf{(a)} 
\begin{align}
 \widetilde{P_n} = P_n, \label{wPntdn}
\end{align}
and
\begin{align}
\widetilde{\beta_n} = \max\{\beta_n, -\ln \ln n\}. \label{wPntdn2}
\end{align}

Given (\ref{wPntdn2}) and $\lim_{n \to \infty}\beta_n = -\infty$, we
clearly obtain $\lim_{n \to \infty}\widetilde{\beta_n} = -\infty$
and $\widetilde{\beta_n} = -O(\ln \ln n)$.

It holds from (\ref{wPntdn2}) that $\widetilde{\beta_n}  \geq
\beta_n$, which along with  (\ref{al0-parta-Hs-od})
 (\ref{al0-parta-Hs}) and (\ref{wPntdn}) yields $t_n \leq \widetilde{t_n}$. Under $t_n \leq \widetilde{t_n}$ and $ \widetilde{P_n} = P_n$, by \cite[Section 3]{zz}, there exists a graph coupling under which $H_s(n,t_n,P_n)$ is a spanning subgraph of $H_s(n,\widetilde{t_n},\widetilde{P_n}) $.

 \noindent
\textbf{(b)} We set
\begin{align}
 \widehat{P_n} = P_n, \label{wPntdn-pb}
\end{align}
and
\begin{align}
\widehat{\beta_n} = \min\{\beta_n, \ln \ln n\}. \label{wPntdn2-pb}
\end{align}

Given (\ref{wPntdn2-pb}) and $\lim_{n \to \infty}\beta_n = \infty$,
we clearly obtain $\lim_{n \to \infty}\widehat{\beta_n} = \infty$
and $\widehat{\beta_n} = O(\ln \ln n)$.

It holds from (\ref{wPntdn2-pb}) that $\widehat{\beta_n}  \leq
\beta_n$, which along with (\ref{al0-parta-Hs-pb-od})
 (\ref{al0-parta-Hs-pb}) and (\ref{wPntdn-pb})  yields $t_n \geq \widehat{t_n}$. Under $t_n \geq \widehat{t_n}$ and $ \widehat{P_n} = P_n$, by \cite[Section 3]{zz}, there exists a graph coupling under which $H_s(n,t_n,P_n)$ is a spanning supergraph of $H_s(n,\widehat{t_n},\widehat{P_n}) $.

\subsection{Proof of Lemma \ref{lem:logn2}.}

From condition
\begin{align}
\frac{1}{s!} \cdot \frac{{K_n}^{2s}}{{P_n}^{s}}   = \frac{\ln  n \pm
O(\ln \ln n) }{n} \sim \frac{\ln  n}{n}, \label{eqslnn-Kpsn}
\end{align}
it holds that
\begin{align}
\frac{{K_n}^2}{P_n}\  & = \Theta\big( n^{-\frac{1}{s}} (\ln
n)^{\frac{1}{s}} \big) , \label{tp-KnPn}
\end{align}
which along with condition $P_n = \Omega(n^c)$ yields
\begin{align}
K_n  &  = \sqrt{P_n \cdot \Theta\big( n^{-\frac{1}{s}} (\ln
n)^{\frac{1}{s}} \big)} = \Omega\Big( n^{\frac{c}{2} - \frac{1}{2s}}
(\ln n)^{\frac{1}{2s}} \Big) . \label{neq-om-Kn2s}
\end{align}

\subsection{Proof of Lemma \ref{qn-dist}.}~

\noindent \textbf{(a)} We still have (\ref{eqslnn-Kpsn}) and (\ref{tp-KnPn}) here. Then
setting $c$ as $1$ in (\ref{neq-om-Kn2s}), it holds that
\begin{align}
K_n  &  = \Omega\Big( n^{\frac{1}{2} - \frac{1}{2s}} (\ln
n)^{\frac{1}{2s}} \Big) . \label{tp-Knsim}
\end{align}
Given (\ref{tp-KnPn}) and (\ref{tp-Knsim}), we use \cite[Lemma
1]{qcomp_kcon} and \cite[Lemma 8]{ZhaoYaganGligor} to have
\begin{align}
q_n    &  = \begin{cases} \frac{1}{s!} \big( \frac{{K_n}^2}{P_n}
\big)^{s}   \big[1\hspace{-1pt}\pm\hspace{-1pt}
O\big(\frac{{K_n}^2}{P_n}\big) \hspace{-1pt}\pm\hspace{-1pt}
O\big(\frac{1}{K_n}\big)\big], &\textrm{\hspace{-2pt}for }s\geq 2,
\vspace{3pt}
\\ \frac{{K_n}^2}{P_n}   \big[1 \hspace{-1pt}\pm\hspace{-1pt}
O\big(\frac{{K_n}^2}{P_n}\big)\big], &\textrm{\hspace{-2pt}for }s=1.
\end{cases} \label{tp-Knsim-qnn}
\end{align}

Now we use (\ref{tp-KnPn}) (\ref{tp-Knsim}) and (\ref{tp-Knsim-qnn})
to derive $q_n  \sim \frac{1}{s!} \cdot
\frac{{K_n}^{2s}}{{P_n}^{s}}$ and $\big| \hspace{2pt} q_n -
\frac{1}{s!} \cdot \frac{{K_n}^{2s}}{{P_n}^{s}} \hspace{2pt} \big|
=   o\big(\frac{1}{n}\big)$.

First, (\ref{tp-KnPn}) and (\ref{tp-Knsim}) imply
$\frac{{K_n}^2}{P_n} = o(1)$ and $K_n = \omega(1)$, respectively,
which are used in (\ref{tp-Knsim-qnn}) to derive $q_n  \sim
\frac{1}{s!} \cdot \frac{{K_n}^{2s}}{{P_n}^{s}}$. Second, applying
(\ref{tp-KnPn}) and (\ref{tp-Knsim}) directly to
(\ref{tp-Knsim-qnn}), we obtain the following two cases:

(i) For $s\geq 2$, it holds that
\begin{align}
&\bigg|q_n -\frac{1}{s!} \bigg( \frac{{K_n}^2}{P_n} \bigg)^{s}\bigg|
\nonumber \\ &  =  \Theta\bigg(\hspace{-1pt}\frac{\ln
n}{n}\hspace{-1pt} \bigg) \Big[\hspace{-2pt}\pm\hspace{-1pt}
O\big(n^{-\frac{1}{s}}(\ln n)^{\frac{1}{s}}\big)
\hspace{-1pt}\pm\hspace{-1pt} O\Big( n^{\frac{1}{2s} - \frac{1}{2}}
(\ln n)^{-\frac{1}{2s}} \hspace{-1pt}\Big)\hspace{-1pt}\Big]
\nonumber \\ & = \pm o\bigg(\frac{1}{n}\bigg)  . \nonumber
\end{align}

(ii) For $s=1$, it holds that
\begin{align}
&\bigg|q_n-\frac{{K_n}^2}{P_n} \bigg|    =  \pm
O\Bigg(\bigg(\frac{\ln n}{n} \bigg)^2 \Bigg) = \pm
o\bigg(\frac{1}{n}\bigg)   . \nonumber
\end{align}

Summarizing cases (i) and (ii) above, we have proved property (a) of
Lemma \ref{qn-dist}.

\noindent \textbf{(b)} By \cite[Lemma 6]{bloznelis2013}, the edge probabilty $q_n$
satisfies
\begin{align}
q_n \leq \frac{\big[\binom{K_n}{s}\big]^2}{\binom{P_n}{s}}.
\label{qnKnPn}
\end{align}
From (\ref{qnKnPn}) and condition $q_n = \frac{\ln  n \pm O(\ln \ln
n) }{n}$, it holds that
\begin{align}
\frac{\big[\binom{K_n}{s}\big]^2}{\binom{P_n}{s}} \geq  \frac{\ln
n}{n} \cdot [1-o(1)] , \nonumber
\end{align}
which along with $\frac{\big[\binom{K_n}{s}\big]^2}{\binom{P_n}{s}}
\leq \frac{1}{s!}  \cdot \frac{{K_n}^{2s} }{(P_n - s +1)^s}$ and
$P_n = \Omega(n)$ leads to
\begin{align}
 {K_n}^{2s}  &  \geq s! (P_n \hspace{-1pt}-\hspace{-1pt} s \hspace{-1pt}+\hspace{-1pt}1)^s \cdot \frac{\ln  n}{n} [1\hspace{-1pt}-\hspace{-1pt}o(1)]  = \Omega(n^{s-1} \ln n) . \nonumber
\end{align}
Therefore, it follows that
\begin{align}
K_n  &  = \Omega\Big( n^{\frac{1}{2} - \frac{1}{2s}} (\ln
n)^{\frac{1}{2s}} \Big) . \label{tp-Knsim-cb}
\end{align}

Note that for some $n$, if $P_n < 2K_n -s $, then two vertices share
at least $s$ items with probability $1$, resulting in $q_n = 1$.
Therefore, given condition $q_n = \frac{\ln  n \pm O(\ln \ln n)
}{n}$, we know that for all $n$ sufficiently large, $P_n \geq 2K_n
-s $ holds, so the probability that two vertices share exactly $s$
items is expressed by ${\binom{K_n}{s} \binom{P_n - K_n}{K_n -
s}}\big/{\binom{P_n}{K_n}}$. Then
\begin{align}
\label{qnnewptv} q_n & \geq  \mathbb{P}[\hspace{2pt}\textrm{Two vertices share exactly $s$ items.}\hspace{2pt}]  \\
& = {\binom{K_n}{s} \binom{P_n - K_n}{K_n - s}}\bigg/{\binom{P_n}{K_n}} \nonumber \\
& = \frac{1}{s!} \cdot \bigg[\prod_{i=0}^{s-1}(K_n-i)\bigg]^2 \cdot
\frac{\prod_{i=0}^{K_n- s-1}(P_n-K_n)}{\prod_{i=0}^{K_n-1}(P_n-i)  }
\nonumber \\
& \geq  \frac{1}{s!} \cdot \frac{(K_n-s+1)^{2s}}{{P_n}^s}  ,
\nonumber
\end{align}
which together with condition $q_n = \frac{\ln  n \pm O(\ln \ln n)
}{n}$ implies
\begin{align}
 \frac{(K_n-s+1)^2}{P_n} = O\big( n^{-\frac{1}{s}} (\ln n)^{\frac{1}{s}} \big).  \label{tp-Knsim-cb2}
\end{align}
From (\ref{tp-Knsim-cb}) and the fact that $s$ is a constant, it
holds that $K_n-s+1 \sim K_n$, which with (\ref{tp-Knsim-cb2})
yields
\begin{align}
 \frac{{K_n}^2}{P_n} = O\big( n^{-\frac{1}{s}} (\ln n)^{\frac{1}{s}} \big).  \label{tp-Knsim-cb3}
\end{align}

Now we use (\ref{tp-Knsim-cb3}) (\ref{tp-Knsim-cb}) and
(\ref{tp-Knsim-qnn}) to derive $q_n  \sim \frac{1}{s!} \cdot
\frac{{K_n}^{2s}}{{P_n}^{s}}$ and $\big| \hspace{2pt} q_n -
\frac{1}{s!} \cdot \frac{{K_n}^{2s}}{{P_n}^{s}} \hspace{2pt} \big|
=   o\big(\frac{1}{n}\big)$, in a way similar to proving property
(a) above.

First, (\ref{tp-Knsim-cb3}) and (\ref{tp-Knsim-cb}) imply $K_n =
\omega(1)$ and $\frac{{K_n}^2}{P_n} = o(1)$, respectively, which are
used in (\ref{tp-Knsim-qnn}) to derive $q_n  \sim \frac{1}{s!} \cdot
\frac{{K_n}^{2s}}{{P_n}^{s}}$. Second, applying (\ref{tp-Knsim-cb3})
and (\ref{tp-Knsim-cb}) directly to (\ref{tp-Knsim-qnn}), we still
have cases (i) and (ii) in the proof of property (a) above. Hence,
finally we also obtain $\big| \hspace{2pt} q_n - \frac{1}{s!} \cdot
\frac{{K_n}^{2s}}{{P_n}^{s}} \hspace{2pt} \big|     =
o\big(\frac{1}{n}\big)$. Then property (b) is proved.

 \subsection{Proof of Lemma \ref{lem_prob_Eij_S1r}.}

 Recall that $E_{ij}$ denotes the event that an edge exists between vertices $v_i$ and $v_j$, and $S_i$ is the number of items on vertex $v_i$. Event $E_{ij}$ occurs if and only if $|S_i \cap S_j| \geq s$. Therefore, event $\bigcap_{j=1}^{r}\overline{E_{ij}} $ is equivalent to $\bigcap_{j=1}^{r} \big(|S_i \cap S_j| < s\big)$, which clearly is implied by event $ \big|S_i \cap \big(\bigcup_{j=1}^{r} S_j\big)\big| < s $.
Then
\begin{align}
\label{Pj1r} &
\mathbb{P}\bigg[\hspace{2pt}\bigcap_{j=1}^{r}\overline{E_{ij}}
\hspace{2pt} \bigg| \hspace{2pt} S_1, S_2, \ldots, S_r\bigg]
  \\ & \leq  \mathbb{P}\bigg[ \hspace{2pt} \bigg|S_i \cap \bigg(\bigcup_{j=1}^{r} S_j\bigg)\bigg| < s \bigg| \hspace{2pt} S_1, S_2, \ldots, S_r \hspace{2pt}\bigg] 
\nonumber \\
& \leq 1 - \mathbb{P}\bigg[ \hspace{2pt}\bigg|S_i \cap
\bigg(\bigcup_{j=1}^{r} S_j\bigg)\bigg| = s \bigg| \hspace{2pt} S_1,
S_2, \ldots, S_r\hspace{2pt}\bigg]
\nonumber \\
& = 1 - \frac{\binom{|\bigcup_{j=1}^{r} S_j|}{s}\binom{P_n-s}{K_n-s}}{\binom{P_n}{K_n}} \nonumber \\
& = 1 - \frac{\binom{|\bigcup_{j=1}^{r} S_j|}{s}\binom{K_n}{s}}{\binom{P_n}{s}} \nonumber \\
& \leq e^{- \frac{\binom{|\bigcup_{j=1}^{r}
S_j|}{s}\binom{K_n}{s}}{\binom{P_n}{s}}} . \nonumber
\end{align}

First, we have (\ref{qnKnPn}) by \cite[Lemma 6]{bloznelis2013}.
Applying (\ref{qnKnPn}) to (\ref{Pj1r}), we obtain
\begin{align}
 \mathbb{P}\bigg[\hspace{2pt}\bigcap_{j=1}^{r}\overline{E_{ij}} \hspace{2pt} \bigg| \hspace{2pt} S_1, S_2, \ldots, S_r\bigg]
& \leq   e^{-\frac{\binom{|\bigcup_{j=1}^{r}
S_j|}{s}}{\binom{K_n}{s}} q_n} . \label{Pj1r2}
\end{align}

Now we prove properties (a), (b), and (c) of Lemma
\ref{lem_prob_Eij_S1r}, respectively.

(a) Given condition $K_n = \omega(1)$, it follows that $\lfloor
(1+{\varepsilon_1}) K_n \rfloor > s $ for all $n$ sufficiently
large. For property (a), we have condition $|\bigcup_{j=1}^{r} S_j|
\geq \lfloor (1+{\varepsilon_1}) K_n \rfloor$, which is used in
(\ref{Pj1r2}) to derive
\begin{align}
 \mathbb{P}\bigg[\hspace{2pt}\bigcap_{j=1}^{r}\overline{E_{ij}} \hspace{2pt} \bigg| \hspace{2pt} S_1, S_2, \ldots, S_r\bigg]
& \leq
e^{-\frac{\binom{\lfloor(1+\varepsilon_1)K_n\rfloor}{s}}{\binom{K_n}{s}}
q_n}. \label{Pj1r3}
\end{align}
We have
\begin{align}
 \label{epsKns} \frac{\binom{\lfloor(1+\varepsilon_1)
K_n\rfloor}{s}}{\binom{K_n}{s}} & =
\frac{\prod_{i=0}^{s-1}\big\{\lfloor(1+\varepsilon_1) K_n\rfloor -
i\big\}}{\prod_{i=0}^{s-1}(K_n-i)}   \\
&  \geq \bigg[\frac{(1+\varepsilon_1) K
_n-1-s}{K_n}\bigg]^s.\nonumber
\end{align}

 Given conditions $ \varepsilon_2 < (1+\varepsilon_1)^s-1$ and $K_n = \omega(1)$, it follows that $K_n \geq  \frac{s+1}{1+\varepsilon_1 - \sqrt[s]{1+\varepsilon_2}} $ for all $n$ sufficiently large, yielding
\begin{align}
\label{epsKns2} \frac{(1+\varepsilon_1) K _n-1-s}{K_n}~~~~~~~~~~~~~~~~~~~~~~~~~~~~~~~~ \\
\geq (1+\varepsilon_1)- (s+1) \cdot \frac{1+\varepsilon_1 -
\sqrt[s]{1+\varepsilon_2}}{s+1} =  \sqrt[s]{1+\varepsilon_2}.
\nonumber
\end{align}
Applying (\ref{epsKns2}) to (\ref{epsKns}), we obtain
\begin{align}
\frac{\binom{\lfloor(1+\varepsilon_1)K_n\rfloor}{s}}{\binom{K_n}{s}}
&    \geq \big(  \sqrt[s]{(1+\varepsilon_2)} \hspace{2pt} \big)^s =
1+\varepsilon_2, \nonumber
\end{align}
which is substituted into (\ref{Pj1r3}) to induce
\begin{align}
 \mathbb{P}\bigg[\hspace{2pt}\bigcap_{j=1}^{r}\overline{E_{ij}} \hspace{2pt} \bigg| \hspace{2pt} S_1, S_2, \ldots, S_r\bigg]
& \leq   e^{- q_n ( 1+\varepsilon_2)}. \nonumber
\end{align}

(b) Given condition $K_n = \omega(1)$, it follows that $\lfloor
\lambda_1 r K_n \rfloor > s $ for all $n$ sufficiently large. For
property (b), we have condition $|\bigcup_{j=1}^{r} S_j| \geq
\lfloor \lambda_1 r K_n \rfloor$, which is used in (\ref{Pj1r2}) to
derive
\begin{align}
 \mathbb{P}\bigg[\hspace{2pt}\bigcap_{j=1}^{r}\overline{E_{ij}} \hspace{2pt} \bigg| \hspace{2pt} S_1, S_2, \ldots, S_r\bigg]
& \leq  e^{-\frac{\binom{\lfloor \lambda_1 r K_n
\rfloor}{s}}{\binom{K_n}{s}} q_n}. \label{Pj1r3b}
\end{align}
We have
\begin{align}
\frac{\binom{\lfloor\lambda_1 r K_n\rfloor}{s}}{\binom{K_n}{s}}
&\hspace{-2pt} =\hspace{-2pt}
\frac{\prod_{i=0}^{s-1}(\lfloor\lambda_1 r K_n \rfloor \hspace{-1pt}
- \hspace{-1pt} i)}{\prod_{i=0}^{s-1}(K_n-i)}  \hspace{-2pt} \geq
\hspace{-2pt} \bigg(\frac{\lambda_1 r K
_n\hspace{-1pt}-\hspace{-1pt}1\hspace{-1pt}-\hspace{-1pt}s}{K_n}\bigg)^s\hspace{-2pt}.
\label{epsKnsb}
\end{align}

 Given conditions $ \lambda_2 < {\lambda_1}^s$ and $K_n = \omega(1)$, it follows that $K_n \geq  \frac{s+1}{\lambda_1 - \sqrt[s]{\lambda_2}} \geq \frac{s+1}{r(\lambda_1 - \sqrt[s]{\lambda_2})}$ for all $n$ sufficiently large, inducing
\begin{align}
\frac{\lambda_1 r K
_n\hspace{-1.5pt}-\hspace{-1.5pt}1\hspace{-1.5pt}-\hspace{-1.5pt}s}{K_n}
\hspace{-1.5pt}\geq\hspace{-1.5pt} \lambda_1 r\hspace{-1.5pt}
-\hspace{-1.5pt} (s\hspace{-1.5pt}+\hspace{-1.5pt}1)
\frac{r(\lambda_1 \hspace{-1.5pt}-\hspace{-1.5pt}
\sqrt[s]{\lambda_2}\hspace{2pt})}{s+1}
\hspace{-1.5pt}=\hspace{-1.5pt} \sqrt[s]{\lambda_2} r.
\label{epsKns2b}
\end{align}
Applying (\ref{epsKns2b}) to (\ref{epsKnsb}), we obtain $ {\binom{\lfloor\lambda_1 r K_n \rfloor}{s}}\Big/{\binom{K_n}{s}} 
\geq  \big(  \sqrt[s]{\lambda_2} r \hspace{2pt} \big)^s =  \lambda_2
r^s \geq \lambda_2 r,$ which is substituted into (\ref{Pj1r3b}) to induce
\begin{align}
 \mathbb{P}\bigg[\hspace{2pt}\bigcap_{j=1}^{r}\overline{E_{ij}} \hspace{2pt} \bigg| \hspace{2pt} S_1, S_2, \ldots, S_r\bigg]
& \leq     e^{-\lambda_2 r q_n}. \nonumber
\end{align}

(c) From $P_n \geq K_n = \omega(1)$, it follows that $P_n =
\omega(1)$. Then
 $\lfloor\mu_1 P_n \rfloor > s $ for all $n$ sufficiently large. For property (c), we have condition $|\bigcup_{j=1}^{r} S_j| \geq  \lfloor\mu_1 P_n\rfloor$, which is used in (\ref{Pj1r}) to derive
\begin{align}
 \mathbb{P}\bigg[\hspace{2pt}\bigcap_{j=1}^{r}\overline{E_{ij}} \hspace{2pt} \bigg| \hspace{2pt} S_1, S_2, \ldots, S_r\bigg]
& \leq  e^{-\frac{\binom{\lfloor\mu_1 P_n
\rfloor}{s}\binom{K_n}{s}}{\binom{P_n}{s}} }. \label{Pj1r3c}
\end{align}
We have
\begin{align}
\label{epsKnsc} & \frac{\binom{\lfloor\mu_1
P_n\rfloor}{s}\binom{K_n}{s}}{\binom{P_n}{s}}   \\ &   \geq
\frac{(s!)^{-1}(\lfloor\mu_1 P_n\rfloor -s )^s \cdot
(s!)^{-1}(K_n-s)^s}{(s!)^{-1}(P_n)^s} \nonumber \\ &   \geq
\frac{1}{s!} \cdot \bigg(\frac{\mu_1 P_n -1 -s}{P_n}\bigg)^s \cdot
(K_n-s)^s.\nonumber
\end{align}
 Given $0 < \mu_2 < (s!)^{-1}{\mu_1}^s$ and $P_n \geq K_n = \omega(1)$, it follows that $P_n \geq \frac{s+1}{\mu_1 - \sqrt{\mu_1}\sqrt[2s]{s! \mu_2}}$ and $K_n \geq \frac{s+1}{1 - \sqrt{\frac{s! \mu_2}{{\mu_1}^s}} }$ for all $n$ sufficiently large, inducing
\begin{align}
  \label{epsKns2c} \frac{\mu_1 P_n - 1 -s}{P_n}~~~~~~~~~~~~~~~~~~~~~~~~~~~~~~~~~~~~~~   \\  \geq \mu_1 - (s+1)
\cdot \frac{\mu_1 - \sqrt{\mu_1}\sqrt[2s]{s! \mu_2}}{s+1} =
\sqrt{\mu_1}\sqrt[2s]{s! \mu_2},\nonumber
\end{align}
and
\begin{align}
 \label{epsKns2c2} (K_n-s)^s~~~~~~~~~~~~~~~~~~~~~~~~~~~~~~~~~~~~~~~~~~~  \\  \geq K_n-s \geq K_n-K_n\Bigg(1 - \sqrt{\frac{s!
\mu_2}{{\mu_1}^s}}\hspace{2pt}\Bigg) = \sqrt{\frac{s!
\mu_2}{{\mu_1}^s}} K_n.\nonumber
\end{align}
Applying (\ref{epsKns2c}) and (\ref{epsKns2c2}) to (\ref{epsKnsc}),
we obtain
\begin{align}
\frac{\binom{\lfloor\mu_1
P_n\rfloor}{s}\binom{K_n}{s}}{\binom{P_n}{s}} \hspace{-1pt}  &
\hspace{-1pt}  \geq \hspace{-1pt}  \frac{1}{s!} \hspace{-.5pt}
 (\sqrt{\mu_1}\sqrt[2s]{s! \mu_2})^s \hspace{-1pt}  \sqrt{\frac{s!
\mu_2}{{\mu_1}^s}} K_n \hspace{-1pt}  = \hspace{-1pt}  \mu_2 K_n,
\nonumber
\end{align}
which is substituted into (\ref{Pj1r3c}) to induce
\begin{align}
 \mathbb{P}\bigg[\hspace{2pt}\bigcap_{j=1}^{r}\overline{E_{ij}} \hspace{2pt} \bigg| \hspace{2pt} S_1, S_2, \ldots, S_r\bigg]
& \leq     e^{-\mu_2 K_n}. \nonumber
\end{align}

 \subsection{Proof of Lemma \ref{olp_lem1}.}

 We consider events $\mathcal{B}_{\ell,r} $, $ \mathcal{D}_{\ell,r}$ and $\mathcal{A}_{{\ell},r}$ defined in Section 6.3.
By definitions, we have
\begin{align}
\mathcal{B}_{\ell,r}&:=  \bigcap_{i=r+1}^{r+\ell} \bigcup_{j =
1}^{r}E_{ij}, \nonumber
\\   \mathcal{D}_{\ell,r}  &:= \bigcap_{i =
r+\ell+1}^{n}\bigcap_{j = 1}^{r} \overline{E_{ij}}, \nonumber
\end{align}
and
\begin{align}
\mathcal{A}_{\ell,r} &:= \mathcal{B}_{\ell,r} \cap \mathcal{C}_{r}
\cap \mathcal{D}_{\ell,r}  . \nonumber
\end{align}
Then considering that given $S_1,S_2, \ldots, S_r$, events $\mathcal{B}_{\ell,r}$ and  $\mathcal{D}_{\ell,r} ~\cap~
\overline{\mathcal{E} (\boldsymbol{J})}$ are conditionally independent, we obtain
\begin{align}
& \label{allbounds} \bP{ \mathcal{A}_{\ell, r} \cap \overline{\mathcal{E} (\boldsymbol{J})}}
 \\
& = \bP{ \mathcal{C}_{r} \cap \mathcal{B}_{\ell,r} \cap \mathcal{D}_{\ell,r} \cap
\overline{\mathcal{E} (\boldsymbol{J})} } \nonumber \\
& = \sum_{\begin{subarray} ~S_1,S_2, \ldots, S_r: \\ ~\mathcal{C}_{r} \textnormal{ happens.} \end{subarray}}\bigg\{ \mathbb{P}[  S_1,\hspace{-1pt} S_2,\hspace{-1pt}
  \ldots,\hspace{-1pt} S_r] \mathbb{P}[ \mathcal{B}_{\ell,r}  \hspace{-1pt} \boldsymbol{\mid}\hspace{-1pt}  S_1,\hspace{-1pt} S_2,\hspace{-1pt}
  \ldots,\hspace{-1pt} S_r] \nonumber \\
&  ~~~~~~~~~~~~~\mathbb{P}[\hspace{2pt}  \mathcal{D}_{\ell,r} ~\cap~
\overline{\mathcal{E} (\boldsymbol{J})}    \hspace{2pt}
\boldsymbol{\mid} \hspace{2pt} S_1, S_2, \ldots, S_r] \bigg\}. \nonumber
\end{align}

We have
 \begin{align}
 \mathbb{P}[ \mathcal{B}_{\ell,r}  \hspace{-1pt} \boldsymbol{\mid}\hspace{-1pt}  S_1,\hspace{-1pt} S_2,\hspace{-1pt}
  \ldots,\hspace{-1pt} S_r]
& \hspace{-1pt} = \hspace{-1pt}  \Bigg\{ \hspace{-1pt}
\mathbb{P}\bigg[ \hspace{-1pt}\bigcup_{j=1}^{r} E_{ij} \hspace{1pt}
\bigg| \hspace{1pt} S_1,\hspace{-1pt} S_2,\hspace{-1pt}
\ldots,\hspace{-1pt} S_r\hspace{-1pt}\bigg]\hspace{-1pt}
\Bigg\}^{\ell}. \nonumber
\end{align}
By the union bound,
\begin{align}
 &\mathbb{P}\bigg[\hspace{2pt}\bigcup_{j=1}^{r} E_{ij} \hspace{2pt} \bigg| \hspace{2pt} S_1, S_2, \ldots, S_r\bigg]
\nonumber \\ &\leq \sum_{j=1}^{r}  \mathbb{P}[ E_{ij}  \hspace{2pt}
\boldsymbol{\mid} \hspace{2pt} S_1, S_2, \ldots, S_r]   =
\sum_{j=1}^{r}  \mathbb{P}[ E_{ij} ]  = r q_n. \nonumber
\end{align}
Then
 \begin{align}
 \mathbb{P}[\hspace{2pt}  \mathcal{B}_{\ell,r}  \hspace{2pt}  \boldsymbol{\mid} \hspace{2pt} S_1, S_2, \ldots, S_r]
& \leq \min\{(r q_n)^{\ell}, 1\} . \label{boundBrleq}
\end{align}
We have
 \begin{align}
 & \mathbb{P}[\hspace{2pt}  \mathcal{D}_{\ell,r} ~\cap~
\overline{\mathcal{E} (\boldsymbol{J})}    \hspace{2pt}
\boldsymbol{\mid} \hspace{2pt} S_1, S_2, \ldots, S_r] \nonumber\\ &
= \Bigg\{ \mathbb{P}\bigg[\hspace{2pt} \bigg(\bigcap_{j=1}^{r}
\overline{E_{ij}} \bigg)~\cap~ \overline{\mathcal{E}
(\boldsymbol{J})}  \hspace{2pt} \bigg| \hspace{2pt} S_1, S_2,
\ldots, S_r\bigg]  \Bigg\}^{n-\ell-r}.\nonumber
\end{align}
By Lemma \ref{lem_prob_Eij_S1r}, for all $n$ sufficiently
large,
\begin{itemize}[leftmargin=20pt]
\item[(a)] for $r = 2, 3, \ldots, R$, it holds that
 \begin{align}
 & \mathbb{P}[\hspace{2pt}  \mathcal{D}_{\ell,r} ~\cap~
\overline{\mathcal{E} (\boldsymbol{J})}    \hspace{2pt}
\boldsymbol{\mid} \hspace{2pt} S_1, S_2, \ldots, S_r] \nonumber \\ &
\quad \leq e^{- q_n (1+\varepsilon_2)(n-\ell-r)} \leq e^{- q_n n
(1+\varepsilon_3)}.\nonumber
\end{align}
To see this, pick any $\varepsilon_3 < (1+ \varepsilon_1)^s-1$, and
use Lemma \ref{lem_prob_Eij_S1r} with $\varepsilon_3 < \varepsilon_2 < (1+ \varepsilon_1)^s-1$.
\item[(b)] for $r = 2, 3, \ldots,  r_n $, it holds that
\begin{align}
&  \mathbb{P}[\hspace{2pt}  \mathcal{D}_{\ell,r} ~\cap~
\overline{\mathcal{E} (\boldsymbol{J})}   \hspace{2pt}
\boldsymbol{\mid} \hspace{2pt} S_1, S_2, \ldots, S_r]  \nonumber \\
& \quad \leq   e^{- \lambda_2 r q_n (n-\ell-r)} \leq e^{- \lambda_2
r q_n n /3}.\nonumber
\end{align}
\item[(c)] for $r = r_n + 1, r_n + 2, \ldots, \lfloor
\frac{n-{\ell}}{2} \rfloor $, it holds that
\begin{align}
&  \mathbb{P}[\hspace{2pt}  \mathcal{D}_{\ell,r} ~\cap~
\overline{\mathcal{E} (\boldsymbol{J})}   \hspace{2pt}
\boldsymbol{\mid} \hspace{2pt} S_1, S_2, \ldots, S_r]  \nonumber \\
& \quad \leq   e^{- \mu_2 K_n (n-\ell-r)} \leq e^{- \mu_2  K_n  n
/3}.\nonumber
\end{align}
\end{itemize}

For simplicity, we use $\Lambda$ to summarize the upper bounds on $\mathbb{P}[\hspace{2pt}  \mathcal{D}_{\ell,r} ~\cap~
\overline{\mathcal{E} (\boldsymbol{J})}    \hspace{2pt}
\boldsymbol{\mid} \hspace{2pt} S_1, S_2, \ldots, S_r]$ in cases (a) (b) and (c) above; i.e., $\Lambda = e^{- q_n n
(1+\varepsilon_3)}$ for $r = 2, 3, \ldots, R$, and $e^{- \lambda_2
r q_n n /3}$ for $r = R+1, R+2, \ldots,  r_n $, and $e^{- \mu_2  K_n  n
/3}$ for $r = r_n + 1, r_n + 2, \ldots, \lfloor
\frac{n-{\ell}}{2} \rfloor $. In view of $\Lambda$, (\ref{boundBrleq}) and $\bP{ \mathcal{C}_{r}}   \leq  \min\{ r^{r-2} {q_n}^{r-1}, 1\}$ by \cite[Lemma 11]{ZhaoYaganGligor}, we obtain from (\ref{allbounds}) that
\begin{align}
& \nonumber \bP{ \mathcal{A}_{\ell, r} \cap \overline{\mathcal{E} (\boldsymbol{J})}}
 \\
& \leq \sum_{\begin{subarray} ~S_1,S_2, \ldots, S_r: \\ ~\mathcal{C}_{r} \textnormal{ happens.} \end{subarray}}\bigg\{ \mathbb{P}[  S_1,\hspace{-1pt} S_2,\hspace{-1pt}
  \ldots,\hspace{-1pt} S_r]\cdot \min\{(r q_n)^{\ell}, 1\} \cdot \Lambda \bigg\} \nonumber \\
&  = \bP{\mathcal{C}_{r}}\cdot  \min\{(r q_n)^{\ell}, 1\} \cdot \Lambda  \\
& \leq  \min\{ r^{r-2} {q_n}^{r-1}, 1\} \cdot  \min\{(r q_n)^{\ell}, 1\} \cdot \Lambda , \nonumber
\end{align}
which clearly completes the proof of Lemma \ref{olp_lem1}.

}

\end{document}